\documentclass{amsart}

\usepackage{amssymb,amsthm,amsmath,amscd}

\def\S{{\mathbb{S}}}

\def\F{{\mathcal{F}}}
\def\eps{\epsilon}
\def\k{\varkappa}
\def\const{\qopname \relax o{const}}
\def\sub{\qopname \relax o{sub}}

\def\an{{\mathcal{A}}_1}

\def\A{{\mathcal{A}}}
\def\ra{\rightarrow}
\def\tr{\qopname \relax o{tr}}

\def\O{\Omega}

\def\E{{\mathbb{E}}}
\def\R{{\mathbb{R}}}
\def\N{{\mathbb{N}}}
\def\Z{{\mathbb{Z}}}

\def\C{{\mathbb{C}}}

\def\Tr{\qopname \relax o{Tr}}
\def\tr{\qopname \relax o{Tr}}
\def\ad{\qopname \relax o{Ad}}

\def\res{\qopname \relax o{Res}}
\def\Res{\qopname \relax o{Res}}
\def\det{\qopname \relax o{det}}
\def\qqquad{\qquad\qquad}

\def\infinity{\infty}
\def\lp{\left(}
\def\rp{\right)}
\def\ls{\left[}
\def\rs{\right]}

\def\vka{{\overline{\varkappa}}}
\def\UpsilonPoio{\Lambda^{(1)}}
\def\UpsilonPoid{\Lambda^{(2)}}
\def\UpsilonPoit{\Lambda^{(3)}}
\def\Upsilonsub{\Upsilon_{3,{{\rm sub}}}}
\def\Upsilonnol{\Upsilon_{3,0}}
\def\UpsilonPoi{\Upsilon_{3,{{\rm Poi}}}}
\def\Upsilondva{\Upsilon_2}
\def\Upsilontri{\Upsilon_3}

\def\spec{\qopname \relax o{spec}}
\def\diag{\qopname \relax o{diag}}
\numberwithin{equation}{section}

\newtheorem{thm}{Theorem}[section]
\newtheorem{lem}[thm]{Lemma}
\newtheorem{cor}[thm]{Corollary}
\newtheorem{prop}[thm]{Proposition}
\theoremstyle{remark}
\newtheorem{rem}{Remark}[section]

\begin{document}

\title{Lower order terms in Szeg\"o type limit 
theorems on Zoll manifolds
}

\author{Dimitri Gioev}
\thanks{The author gratefully acknowledges a full support for the academic year 2002--03 from the Swedish Foundation for International Cooperation in Research and Higher Education (STINT), Dnr.~PD2001--128.}

\address{Department of Mathematics, University of Pennsylvania,
DRL 209 South 33rd Street, Philadelphia, PA~19104--6395}
\curraddr{Department of Mathematics, Courant Institute of Mathematical Sciences,
New York University, 251 Mercer Street, New York, NY~10012--1185}
%\email{gioev@math.upenn.edu}
\email{gioev@math.upenn.edu, gioev@cims.nyu.edu}

\begin{abstract}
We compute the third order term in a
generalization of the Strong Szeg\"o Limit Theorem
for a zeroth order pseudodifferential operator (PsDO) on a Zoll manifold
of an arbitrary dimension. 
In \cite{GO}, the second order term was computed 
by V.~Guillemin and K.~Okikiolu.
In the present paper, 
an important role 
is played by a certain combinatorial identity which we call the generalized
Hunt--Dyson formula \cite{GLc}.
%identity which generalizes the formula of 
%G.~A.~Hunt and F.~J.~Dyson to an arbitrary natural power.
%It turns out that t
This identity is a different form of 
the renowned Bohnenblust--Spitzer combinatorial theorem
which is related to the maximum of a random walk with i.i.d. steps on the real line.
% of H.~F.~Bohnenblust
%that has been employed by F.~Spitzer to compute the characteristic
%function of the maximum of a random walk with independent identically
%distributed steps.
%
A corollary of our main result is a fourth order Szeg\"o type
asymptotics for a zeroth order PsDO on the unit circle, which
in matrix terms gives
a fourth order asymptotic formula
for the determinant of the truncated 
sum $P_n(T_1+T_2D)P_n$
of a Toeplitz matrix $T_1$ with the product 
of another Toeplitz matrix $T_2$ 
and a diagonal matrix $D$
of the form $\diag(0,1,\frac12,\frac13,\cdots)$. 
Here $P_n=\diag(1,\cdots,1,0,\cdots)$, $n$ ones.
\end{abstract}

%\begin{keyword}
%Toeplitz operators \sep Szeg\"o type asymptotics \sep Zoll manifolds
%\end{keyword}
%\dedication{Dedicated to...}
\maketitle
%\section{Introduction and main results}
\section{Introduction}
%\label{CH1_s5intro}
%
The main motivation for this work was to find an
explicit formula for a ``Szeg\"o--regularized'' determinant
of a zeroth order pseudodifferential operator (PsDO)
on a Zoll manifold introduced in \cite[after~(3)]{GO3} and \cite{O2}, 
see Remark~\ref{remdet}.
Our main result, Theorem~\ref{thm_main}, is valid for any dimension
$d\in\N$. In the case $d=2$, Theorem~\ref{thm_main}
gives such a formula.
\subsection{Notations and main results}
Let $M=\S^1$ be the unit circle $\R/2\pi\Z$. 
Denote by $P_n$, $n\in\N$, the orthogonal
projection from $L^2(\S^1)$ 
to the subspace spanned by $\{e^{ikx}\}_{|k|\leq{}n}$. 
For a function $f\in{}L^1(\S^1)$ 
denote its $k$th
Fourier coefficient by
$\widehat{f}_k := \int_0^{2\pi} f(x) e^{-ikx}\,\frac{dx}{2\pi}$, $k\in\Z$.
Let $b(x)$ be a positive function on $\S^1$
such that $\sum_{k\in\Z}|k|\,|\widehat{(\log{b})}_k|^2<\infty$.
Denote by $B$ the operator of multiplication by
$b$ acting in $L^2(\S^1)$.
The matrix representation of the operator $B$
in the basis $\{e^{ikx}\}_{k\in\Z}$ 
is the Toeplitz matrix $(\widehat{b}_{j-k})_{j,k\in\Z}$.
The classical Strong Szeg\"o Limit Theorem (SSLT)~\cite{Sz2}
states that
$$
%\begin{aligned}
   \tr\log P_nBP_n = \tr P_n(\log{B})P_n
              +\sum_{k=1}^\infty k\,\widehat{(\log b)}_k\,
                       \widehat{(\log b)}_{-k}
              +o(1),\qquad n\ra\infty.
%\end{aligned}
$$
Here $\tr\log P_nBP_n=\log\det P_nBP_n$
and $\tr P_n(\log{B})P_n=(2n+1)\int_0^{2\pi} \log b(x)\,\frac{dx}{2\pi}$.
It has been shown by H.~Widom that the remainder
is $O(n^{-\infty})$ if $b(x)\in{}C^{\infty}(\S^1)$, see \cite{Wajm}.

The main result of this paper is Theorem~\ref{thm_main},
in which we find a third order
generalization of the SSLT for a zeroth order 
pseudodifferential operator (PsDO) $B$ on
a Zoll manifold $M$ of an arbitrary dimension
$d\in\N$.

Recall that $M$ is called a {\em Zoll manifold\ }\cite{GO3}
if it is compact, closed
and such that the geodesic flow
on $M$ is simply periodic with period $2\pi$. 
The unit circle and the standard sphere of any dimension
are Zoll manifolds.
A second order generalization of the SSLT for a Zoll manifold $M$
of any dimension has been obtained by V.~Guillemin
and K.~Okikiolu \cite{GO3,GO}, 
see also an important preceding work \cite{O} by K.~Okikiolu
for $M=\S^2$ and $\S^3$.
The proofs in \cite{O,GO3,GO} use a combinatorial
identity due to G.~A.~Hunt and F.~J.~Dyson and proceed in
the spirit of the combinatorial proof of the classical SSLT 
by M.~Kac \cite{K}. See also \cite{GO2,O2} 
where the combinatorial approach and the Hunt--Dyson
formula (HD) are used in a different setting
to obtain
a second order generalization of the SSLT 
for a manifold with the set of closed geodesics
of measure zero in the unit cotangent bundle. 

In the proof of Theorem~\ref{thm_main} we use the method of \cite{GO}.
A central role in our proof is
played by a certain combinatorial identity
which generalizes the Hunt--Dyson 
formula mentioned above to an arbitrary natural power.
We call this identity the {\em generalized Hunt--Dyson formula (gHD)}, 
see Theorem~\ref{CH1_S2TH1} and \cite{GLc}.
After having discovered and proved the gHD
we realized that it is related to 
another combinatorial theorem, which has a long history. 
This theorem 
is a result due to H.~F.~Bohnenblust
that appeared in an article by F.~Spitzer
on random walks \cite[Theorem~2.2]{S1}, 
and is now commonly known as the
{\em Bohnenblust--Spitzer theorem (BSt)}.
A major application of, and motivation for the BSt,
is the computation of 
the characteristic function of the maximum of a random walk
with independent identically distributed (i.i.d.)~steps carried out in \cite{S1}.
Note that the expectation
of such a maximum was computed earlier in \cite{K} with the
help of the usual HD.

Let $M$ be a Zoll manifold of dimension $d\in\N$. 
Let $\Psi^m(M)$, $m\in\Z$, denote
the space of classical PsDO's of order $m$ on $M$. 
Recall that for a given $G\in\Psi^m(M)$,
its principal symbol $\sigma_m(G)$ and subprincipal symbol $\sub(G)$
are well-defined on $T^*M$.
Let $\Delta$ denote the Laplace--Beltrami operator on $M$.
It is known \cite{DG} that there exists a constant $\alpha\in\R$
such that the spectrum of  $\sqrt{-\Delta}$
lies in bands around the points $k+\frac{\alpha}{4}$, $k\in\N$.
Moreover, it has been shown in \cite{CdV}
that there exists $A_{-1}\in\Psi^{-1}(M)$
such that $[\Delta,A_{-1}]=0$
and the spectrum of the operator 
\begin{equation}
\label{eq_A}
     A:= \sqrt{-\Delta}-\frac{\alpha}4 - A_{-1}
\end{equation}
is $\N$.
Let $P_n$, $n\in\N$, denote the projection from $L^2(M)$
onto the subspace spanned by the eigenfunctions
of $A$ corresponding to the eigenvalues $1,2,\cdots,n$. 
%Denote by $dxd\xi$ the standard measure
%on 
Let $dxd\xi$ be the standard measure on 
$S^*M:=\{(x,\xi):\sigma_1(A)(x,\xi)=1\}$ divided by $(2\pi)^d$. 
Following \cite{GO} we will assume 
that $\sigma_1(A)(x,\xi)=\sigma_1(A)(x,-\xi)$
for all $(x,\xi)\in{}T^*M$. In \cite[Chapter~1]{GiPhD} this is not
assumed which leads to more complicated expressions.
Let $\Theta^t(x,\xi)$ denote the shift of the point $(x,\xi)\in{}S^*M$
by $t$ units along geodesic flow.
For any function $f\in{}C^\infty(S^*M)$ introduce
the $k$th Fourier coefficient along the closed geodesic
of length $2\pi$ starting at a given point $(x,\xi)$
\begin{equation}
\label{eq_FcfZ} 
        \widehat{f}_k(x,\xi):= \int_0^{2\pi} e^{-ikt} f(\Theta^t(x,\xi))\,
     \frac{dt}{2\pi},\qquad k\in\Z.
\end{equation}

The simplest form of our result is for the case of $M=\S^1$ 
with $f(z)=\log{}z$. Note that $f(z)$ is
analytic in a disk of radius $1$ about the point $z=1$.
In the proofs in Section~\ref{CH1_s95},
we require that the function $f(z)$ is analytic
on a disk the radius of which depends on
a certain %symbolic 
norm of the operator $B\in\Psi^0(M)$.
For our purposes the following norm is convenient
\begin{equation}
\label{eq_seminorm}
\begin{aligned}
% |||B|||_d := \|B\|+\big\|[A,[A,B]]\big\|
%&+\|\nabla^{\max(d+2,6)}\sigma_0(B)\|_\infty\\
%   &+ \|\nabla^3\sub(B)\|_\infty
%   + \|\nabla^3\sigma_{-2}(B)\|_\infty,
  |||B|||_d := \|B\|&+\big\|[A,[A,B]]\big\|
+\|\nabla^{\max(d+2,6)}\sigma_0(B)\|_\infty\\
   &+ \|\nabla^3\sub(B)\|_\infty
     + \int_{S^*M}\big|\nabla^3\sigma_{-d}(A^{2-d}B)\big|\,dxd\xi,
%   + \|\nabla^3\sigma_{-2}(B)\|_\infty,
\end{aligned}
    %+ \|\nabla^3\sigma_{-d}(A^{2-d}B)\|_\infty,
%\max\big\{\|B\|, \big\|[A,[A,B]]\big\|
%\|\nabla^{\max(d+2,6)}\sigma_0(B)\|_\infty,
%    \|\nabla^3\sub(B)\|_\infty,
%    \|\nabla^3\sigma_{-2}(B)\|_\infty\big\},
%    %+ \|\nabla^3\sigma_{-d}(A^{2-d}B)\|_\infty,
\end{equation}
where $d=\dim{}M$,
$\|\cdot\|$ is the operator norm in $L^2(M)$, $\nabla$
includes both $x$- and $\xi$-derivatives in local coordinates on $T^*M$,
and $\|g\|_\infty:=\max_{(x,\xi)\in{}S^*M}|g(x,\xi)|$. (The integral
in \eqref{eq_seminorm} is well-defined being a Guillemin--Wodzicki
residue, see \eqref{eqR} below.)
%It would be interesting to know what 
\begin{thm}
\label{thm_onedim}
 Let $M=\S^1$ and $P_n$ be the projection
on the linear span of $\{e^{ikx}\}_{|k|\leq{}n}$.
Let $B\in\Psi^0(M)$ and assume that $\sigma_0(B)$
is strictly positive, and that the symbolic norm $|||I-B|||_1$
is sufficiently small. Then $\log{}B\in\Psi^0(M)$ and 
the following holds as $n\ra\infty$,
\begin{equation}
\label{eq_f}
\begin{aligned}
    \tr\log P_nBP_n &= \tr P_n(\log B)P_n 
                + \frac12\int_{S^*M}\sum_{k=1}^\infty k\,\widehat{(\sigma_0(\log{B}))}_k
                   \,\widehat{(\sigma_0(\log{B}))}_{-k}\,dxd\xi\\
              &+\frac1n\cdot
                     \frac12\int_{S^*M}\sum_{k=1}^\infty k\,\widehat{(\sigma_0(\log{B}))}_k
                   \,\widehat{(\sub(\log B))}_{-k}\,dxd\xi
         +O\bigg(\frac1{n^2}\bigg).
\end{aligned}
\end{equation}
\end{thm}
In \eqref{eq_f} the argument $(x,\xi)\in{}S^*M$ is omitted 
for brevity and for each $(x,\xi)\in{}S^*M$ the Fourier coefficient 
is understood in the sense of \eqref{eq_FcfZ}.

We need to fix %some 
more notation to formulate the result for $M=\S^1$
and an arbitrary analytic $f(z)$.
Let $\an$ denote the set of all analytic functions on $\C$ with no
constant term
$$
             \an:=\Big\{f(z)\,:\,f(z)=\sum_{m=1}^\infty c_m z^m,\,z\in\C\Big\}.
$$
%By a linear functor defined on $\an$ w
%We will call a linear map defined$\an$ to the space of 
%For any $f\in\an$ define \cite{LRS}
In \cite{LRS} the authors introduce a linear map $W_2$ 
from $\an$ to the space of continuous functions from $\C^2$ to $\C$, 
defined by
\begin{equation}
\label{CH1_s5ieW}
   W_2[f](x_1,x_2):=\frac12\int_0^{x_1}\!\!\int_0^{x_2}
        \frac{f^\prime(\xi_1)-f^\prime(\xi_2)}{\xi_1-\xi_2}\,d\xi_1 d\xi_2.
\end{equation}
Let $j\in\N$. We will call a linear map $V$ from $\an$ to the space 
of continuous functions $\C^j\rightarrow\C$ a {\em linear $j$-map.\ }The 
linearity means %of course 
that for arbitrary $f,g\in\an$, $\alpha,\beta\in\C$,
$x_1,\cdots,x_j\in\C$
$$
\begin{aligned}
       V[\alpha{}f+\beta{}g](x_1,\cdots,x_j) &= 
         \alpha{}V[f](x_1,\cdots,x_j)%\\
              %&\quad
+\beta{}V[g](x_1,\cdots,x_j).
\end{aligned}
$$
A $2$-map $U$, which is equivalent to $W_2$, was earlier
constructed by H.~Widom \cite{W4,W5}.
We will need also a $2$-map $\tilde{W}_2$ whose action on an
arbitrary $f\in\an$ is prescribed by
\begin{equation}
\label{CH1_s5ieWt}
   \tilde{W}_2[f](x_1,x_2):=\frac12\int_0^{x_1}
        \frac{f^\prime(\xi_1)-f^\prime(x_2)}{\xi_1-x_2}\,d\xi_1.
\end{equation}
For an arbitrary $B\in\Psi^0(M)$, let us write 
$b_0:=\sigma_0(B)$ and $b_{{\rm sub}}:=\sub(B)$,
$b_0^t(x,\xi):=b_0(\Theta^t(x,\xi))$,
$b_{{\rm sub}}^t(x,\xi):=b_{{\rm sub}}(\Theta^t(x,\xi))$,
and omit the argument $(x,\xi)\in{}S^*M$. % for brevity.
 It is convenient
to introduce the notations 
$$%\begin{equation}
%\label{eqszdva}
   \Upsilondva[f](B):=\int_{S^*M}dxd\xi 
     \sum_{k=1}^\infty k\,\int_0^{2\pi}\int_0^{2\pi}
                    e^{ik(t_1-t_2)}\,W_2[f](b_0^{t_1},b_0^{t_2})\,
                  \frac{dt_1}{2\pi}\frac{dt_2}{2\pi}\\
$$%\end{equation}
and
$$%\begin{equation}
%\label{eqsztrisub}
      \Upsilonsub[f](B)  :=
        \int_{S^*M}\,dxd\xi\sum_{k=1}^\infty k\,\int_0^{2\pi}\int_0^{2\pi}
                    e^{ik(t_1-t_2)}\,\tilde{W}_2[f](b_0^{t_1},b_0^{t_2})\,
                  b_{{\rm sub}}^{t_2}\,\frac{dt_1}{2\pi}\frac{dt_2}{2\pi}.
$$%\end{equation}
\begin{thm}
\label{thm_onedimany}
 Let $M=\S^1$ and $P_n$ be the projection
on the linear span of $\{e^{ikx}\}_{|k|\leq{}n}$.
Let $B\in\Psi^0(M)$ and $f\in\an$.
Then $f(B)\in\Psi^0(M)$ and 
the following holds as $n\ra\infty$,
$$%\begin{equation}
%\label{eq_fany}
%\begin{aligned}
    \tr f(P_nBP_n) = \tr P_n f(B)P_n 
              + \Upsilondva[f](B)+\frac1n\cdot\Upsilonsub[f](B)
         +O\bigg(\frac1{n^2}\bigg). %,\quad n\ra\infty.
%\end{aligned}
$$%\end{equation}
\end{thm}
Observe that Theorem \ref{thm_onedim} follows
from Theorem \ref{thm_onedimany} setting $f(t)=\log{}z$ 
and noting following H.~Widom and \cite{LRS,O2} that
\begin{equation}
\label{eq_fact}
   W_2[\log](x_1,x_2)= -\frac12\log x_1 \log x_2,
\end{equation}
and also that for $0\leq{}t_1,t_2\leq2\pi$
$$
      \tilde{W}_2[\log](b_0^{t_1},b_0^{t_2})b_{{\rm sub}}^{t_2} 
         =(\log b_0^{t_1}) (b_{{\rm sub}}^{t_2}/b_0^{t_2})=
          (\sigma_0(\log B))^{t_1}\,(\sub(\log B))^{t_2}.
$$

In the higher dimensional case two additional contributions
to the third Szeg\"o term now arise.
First, define a $3$-map ${W}_3$ such that for any $f\in\an$
\begin{equation}
\label{CH1_s5ieFprin}
\begin{aligned}
     W_3[f](x_1,x_2,x_3)&:=x_3\int_0^{x_1}\int_0^{x_2}
          \bigg(\frac{f(\xi_1)}{\xi_1(\xi_1-x_3)(\xi_1-\xi_2)}\\
     -&\frac{f(\xi_2)}{\xi_2(\xi_2-x_3)(\xi_1-\xi_2)}
           +\frac{f(x_3)}{x_3(\xi_1-x_3)(\xi_2-x_3)}
      \bigg) \,d\xi_1d\xi_2,
\end{aligned}
\end{equation}
and introduce the notation
$$%\begin{equation}
%\label{eqsztrinol}
\begin{aligned}
      \null\Upsilonnol&[f](B)  :=(d-1)\int_{S^*M}dxd\xi\\
             &\times\bigg[
           \sum_{k=1}^\infty\Big(k^2+(1+\alpha/2) k\Big)\,
      \int_0^{2\pi}\int_0^{2\pi}
                    e^{ik(t_1-t_2)}\,W_2[f](b_0^{t_1},b_0^{t_2})\,
                  \frac{dt_1}{2\pi}\frac{dt_2}{2\pi}\\
               &\quad+\sum_{k,l=1}^\infty  k\,l\,\int_0^{2\pi}\!\!
\int_0^{2\pi}\!\!\int_0^{2\pi} e^{ik(t_1-t_3)+il(t_2-t_3)}
    \,W_3[f](b_0^{t_1},b_0^{t_2},b_0^{t_3})\,\frac{dt_1}{2\pi}
\frac{dt_2}{2\pi}\frac{dt_3}{2\pi}\bigg],
\end{aligned}
$$%\end{equation}
see Remark \ref{rem_RR}.
Secondly,  we introduce notations needed to describe the contribution
of the Poisson brackets $\{b_0^t,b_0^s\}$, $0\leq{}t,s\leq2\pi$, 
to the third Szeg\"o term.
Define for each $j\in\N$ a linear $j$-map $\Phi_j$
such that for any $m\in\N$, $m\geq{}j$,
\begin{equation}
\label{CH1_s5iePhi2}
       \Phi_j[z^m](x_1,\cdots,x_j) := \sum_{ \genfrac{}{}{0pt}{}{l_1,\cdots,l_j\geq1}
         {l_1+\cdots+l_j=m} } 
               \frac{x_1^{l_1}}{l_1}\cdots \frac{x_j^{l_j}}{l_j},
\end{equation}
and $\Phi_j[z^m]:=0$ for $m=1,\cdots,j-1$.
This together with the linearity defines $\Phi_j$ uniquely on $\an$, see
Section~\ref{CH1_s6Phi} and  \eqref{CH1_s6Phi.e2p}
for an expression for $\Phi_j$ acting on an arbitrary $f\in\an$.
%Note that all the powers in \eqref{CH1_s5iePhi2} are positive.
We can write \eqref{CH1_s5iePhi2} in terms of the complete symmetric function
of degree $m-j$ evaluated at the point $(\xi_1,{}\cdots,{}\xi_j,{}0,{}\cdots)$,
see \eqref{CH1_s6Phi.e1}.
%%In Section~\ref{CH1_s6Phi}, we give
%%an %invariant (not depending on $m$ as in \eqref{CH1_s5iePhi2})
%%expression for $\Phi_j$ acting on an arbitrary $f\in\an$, see \eqref{CH1_s6Phi.e2p}. 
%
Now for $x,y\in\R$ denote $-(x)_-:=\min(0,x)$ 
and $M_2(x,y):=\min(0,x,x+y)$.
For arbitrary $j,k,l\in\N$ and any $\k_1,\k_2,
           \mu_1,\cdots,\mu_j$, $
          \nu_1,\cdots,\nu_k$, $
         \rho_1,\cdots,\rho_l\in\Z$ 
set
%\begin{equation}
\begin{align}
\label{CH1_s5ieOj1}
%\begin{aligned}
     \null\Omega_j^{(1)}&(\k_1,\k_2,\mu_1,\cdots,\mu_j)\notag\\ 
          &:= M_2(\k_1,\k_2) - \big(\,\k_1+\k_2 - M_2(\k_1,\k_2)
            -(\mu_1)_- -\cdots -(\mu_j)_- \,\big)_-\notag\\
     &-(\k_1)_- -\Big(\,(\k_1)_+-(\mu_1)_--\cdots-(\mu_j)_- \\
                   &-\big(\,(\mu_1)_++\cdots+(\mu_j)_++\k_2\,\big)_-\,\Big)_-\notag\\
   &-(\mu_1)_- -\cdots -(\mu_j)_- 
        -\big(\,(\mu_1)_++\cdots+(\mu_j)_++M_2(\k_1,\k_2)\,\big)_-\notag,
%\end{aligned}
%\end{equation}
\end{align}
and
\begin{align}
%\begin{equation}
\label{CH1_s5ieOj2}
      \null\Omega_{j,k}^{(2)}&(\k_1,\k_2,\mu_1,\cdots,\mu_j,\nu_1,\cdots,\nu_k)\notag\\ 
          &:= -(\k_1)_- -\Big(\,(\k_1)_+-(\mu_1)_--\cdots-(\mu_j)_-\notag\\
                         &-\big(\,(\mu_1)_++\cdots+(\mu_j)_+ +\k_2 
           -(\nu_1)_--\cdots-(\nu_k)_-\,\big)_-\,\Big)_-\notag\\
   &-(\mu_1)_- -\cdots -(\mu_j)_- 
        -\Big(\,(\mu_1)_++\cdots+(\mu_j)_+ +\k_1-(\k_2)_- \notag\\
                                  &-\big(\,(\k_2)_+ 
     -(\nu_1)_--\cdots-(\nu_k)_-\,\big)_-\,\Big)_-\\
   &-(\mu_1)_- -\cdots -(\mu_j)_- 
        -\Big(\,(\mu_1)_++\cdots+(\mu_j)_+ +\k_1  \notag\\
                                     &-(\nu_1)_--\cdots-(\nu_k)_-
          -\big(\,(\nu_1)_++\cdots+(\nu_k)_++\k_2\,\big)_-\,\Big)_-,\notag
\end{align}
%\end{equation}
and
\begin{align}
\label{CH1_s5ieOj3}
     \null\Omega_{j,k,l}^{(3)}&(\k_1,\k_2,\mu_1,\cdots,\mu_j,\nu_1,\cdots,\nu_k,
               \rho_1,\cdots,\rho_l)\notag\\ 
          &:= -(\mu_1)_- -\cdots -(\mu_j)_- 
                 -\Big(\,(\mu_1)_++\cdots+(\mu_j)_+\notag \\
           &+\k_1 -(\nu_1)_--\cdots-(\nu_k)_- 
          - \big(\,(\nu_1)_++\cdots+(\nu_k)_+\\
          &+ \k_2  -(\rho_1)_--\cdots-(\rho_l)_-\,\big)_-\,\Big)_-.\notag
\end{align}
For any $f\in\an$ denote by 
$T_j[f](z)$ its Taylor polynomial of degree $j\in\N$ about the point $t=0$.
For any $f\in\an$ and $B\in\Psi^0(M)$ %, $b_0:=\sigma_0(B)$. 
introduce the notation
$$%\begin{equation}
%\label{eqsztripoi}
      \UpsilonPoi[f](B)  := \int_{S^*M}\Big(\UpsilonPoio[f](b_0)
      + \UpsilonPoid[f](b_0)
       + \UpsilonPoit[f](b_0)\Big)\,dxd\xi
$$%\end{equation}
where 
\begin{equation}
\label{CH1_s5ie1Phi20Poi1}
\begin{aligned}
     \UpsilonPoio[f](b_0)&:=\frac1{2i}\,
      \sum_{j=1}^{\infty}\frac1{j!}\sum_{{\k_1+\k_2+\mu_1+\cdots+\mu_j=0}}
   \Omega_j^{(1)}(\k_1,\k_2,\mu_1,\cdots,\mu_j)\\
&\times\int_0^{2\pi}\!\!\int_0^{2\pi} e^{i(\k_1u_1+\k_2u_2)}\,
\big\{\,b_0^{u_1},b_0^{u_2}\big\}\frac{du_1}{2\pi}\frac{du_2}{2\pi}\\
&\times\int_0^{2\pi}\!\cdots\!\int_0^{2\pi}\,
e^{i(\mu_1r_1+\cdots+\mu_jr_j)}\\
&\times\Phi_j\big[z^{-2}(f(z)-T_2[f](z))\big](b_0^{r_1},\cdots,b_0^{r_j})
\,\frac{dr_1}{2\pi}\cdots\frac{dr_j}{2\pi},
\end{aligned}
\end{equation}
%
%the second is
and
%\begin{equation}
\begin{align}
\label{CH1_s5ie1Phi20Poi2}
%\begin{aligned}
   \UpsilonPoid&[f](b_0):=\frac1{2i}\,\sum_{j,k=1}^{\infty}\frac1{j!}\frac1{k!}\notag\\
      &\times
   \sum_{ \genfrac{}{}{0pt}{}{\k_1+\k_2+\mu_1+\cdots+\mu_j+}
        {+\nu_1+\cdots+\nu_k=0} }
   \Omega_{j,k}^{(2)}(\k_1,\k_2,\mu_1,\cdots,\mu_j,\nu_1,\cdots,\mu_k)
  \notag\displaybreak[0]\\
&\times\int_0^{2\pi}\!\!\int_0^{2\pi} e^{i(\k_1u_1+\k_2u_2)}\,
\big\{\,b_0^{u_1},b_0^{u_2}\big\}\frac{du_1}{2\pi}\frac{du_2}{2\pi}\\
&\times\int_0^{2\pi}\!\cdots\!\int_0^{2\pi}\,
e^{i(\mu_1r_1+\cdots+\mu_jr_j+\nu_1s_1+\cdots+\nu_ks_k)}
\,\frac{dr_1}{2\pi}\cdots\frac{dr_j}{2\pi}
\frac{ds_1}{2\pi}\cdots\frac{ds_k}{2\pi}\notag\\
&\times\Phi_{j+k}\big[z^{-2}(f(z)-T_3[f](z))\big](b_0^{r_1},\cdots,b_0^{r_j},
       b_0^{s_1},\cdots,b_0^{s_k}),\notag
\end{align}
%\end{aligned}
%\end{equation}
%
and 
%\begin{equation}
\begin{align}
%%\label{CH1_s5ie1Phi20Poi3}
%\begin{aligned}
   &\UpsilonPoit[f](b_0):=\frac1{2i}
              \,\sum_{j,k,l=1}^{\infty}
       \frac1{j!}\,\frac1{k!}\,\frac1{l!}\notag\\
   &\,\,\times\sum_{ \genfrac{}{}{0pt}{}{\k_1+\k_2+\mu_1+\cdots+\mu_j+}
    {+\nu_1+\cdots+\nu_k+\rho_1+\cdots+\rho_l=0} }
   \Omega_{j,k,l}^{(3)}(\k_1,\k_2,\mu_1,\cdots,\mu_j,\nu_1,\cdots,\nu_k,
        \rho_1,\cdots,\rho_l)\notag\\
&\,\,\times\int_0^{2\pi}\!\!\int_0^{2\pi} e^{i(\k_1u_1+\k_2u_2)}\,
\label{CH1_s5ie1Phi20Poi3}
\big\{\,b_0^{u_1},b_0^{u_2}\big\}\frac{du_1}{2\pi}\frac{du_2}{2\pi}\\
&\,\,\times\int_0^{2\pi}\!\cdots\!\int_0^{2\pi}
\,\frac{dr_1}{2\pi}\cdots\frac{dr_j}{2\pi}
\,\frac{ds_1}{2\pi}\cdots\frac{ds_k}{2\pi}
\,\frac{dt_1}{2\pi}\cdots\frac{dt_l}{2\pi}\notag\\
&\,\,\times e^{i(\mu_1r_1+\cdots+\mu_jr_j+\nu_1s_1+\cdots+\nu_ks_k
       +\rho_1t_1+\cdots+\rho_lt_l)}\notag\\
&\,\,\times\Phi_{j+k+l}\Big[z^{-2}\big(f(z)-T_4[f](z)\big)\Big]
       (b_0^{r_1},\cdots,b_0^{r_j},
       b_0^{s_1},\cdots,b_0^{s_k},b_0^{t_1},\cdots,b_0^{t_k}).\notag
\end{align}
%\end{aligned}
%\end{equation}
Now the result for any dimension and an arbitrary $f\in\an$. 
Denote
$$%\begin{equation}
%\label{eqsztri}
         \Upsilontri[f](B):=  \Upsilonnol[f](B) + \Upsilonsub[f](B)  
       +     \UpsilonPoi[f](B).
$$%\end{equation}
%\begin{rem}
%An analogous result holds for an arbitrary $f\in\an$,
%simply replace in \eqref{eq_th2} the expressions
%$\widehat{(\log b_0)}_k\,\widehat{(\log b_0)}_{-k}$
%and $\widehat{(\log b_0)}_k\,\widehat{(b_{{\rm sub}}/b_0)}_{-k}$
%by their counterparts from Theorem~\ref{thm_onedimany}.
%\end{rem}
\begin{thm}
\label{thm_main}
Let $M$ be a Zoll manifold of dimension $d\in\N$. 
Let $A$ be defined by \eqref{eq_A}.
Assume that $\sigma_1(A)(x,\xi)=\sigma_1(A)(x,-\xi)$ for all $(x,\xi)\in{}T^*M$.
Let $B\in\Psi^0(M)$ and $f\in\an$. Then $f(B)\in\Psi^0(M)$ and 
the following holds as $n\ra\infty$,
\begin{equation}
\label{eq_th5}
\begin{aligned}
    \tr f(P_nBP_n) = \tr P_n f(B)P_n  &+ n^{d-1}\cdot\Upsilondva[f](B)\\
           &+n^{d-2}\cdot\Upsilontri[f](B)
            +O(n^{d-3}).
\end{aligned}
\end{equation}
\end{thm}
Most of the paper is devoted to the proof of Theorem~\ref{thm_main}.
Theorem~\ref{thm_onedimany} follows 
from Theorem \ref{thm_main} in view of the following. 
For $d=1$, $\Upsilonnol$ vanishes,
and also $\Upsilon_{3,{{\rm Poi}}}$ vanishes, because all the 
Poisson brackets
vanish in this case (for each of the two cotangent directions
the angle does not change and
$\sigma_0(B)$ is homogeneous of degree $0$ in $\xi$).
%(see \cite{W2})
%%Theorem~\ref{thm_onedim} follows from 
%%Theorem~\ref{thm_onedimany}
%%because of \eqref{eq_fact} and the fact that for any $0\leq{}t_1,t_2\leq2\pi$
%%$$
%%      \tilde{W}_2[\log](b_0^{t_1},b_0^{t_2})b_{{\rm sub}}^{t_2} 
%%         =(\log b_0^{t_1}) (b_{{\rm sub}}^{t_2}/b_0^{t_2})=
%%          (\sigma_0(\log B))^{t_1}\,(\sub(\log B))^{t_2}.
%%$$
Note that for the function $f(z)=\log{}z$,
the terms involving $W_2$ and $\tilde{W}_2$
in Theorem~\ref{thm_main} for any $d\in\N$
take a simpler form as in Theorem~\ref{thm_onedim}.
We make several remarks, and then state a corollary to 
Theorem \ref{thm_onedim} and \ref{thm_main}
which gives an explicit formula for $\log\det P_nBP_n$,
as $n\ra\infty$. For $d=1,2$, that formula gives an expression
for a possible generalized determinant of the operator $B\in\Psi^0(M)$,
see Remark \ref{remdet} below. 
\begin{rem}
The {\em existence\ }of a full expansion of the type \eqref{eq_th5}
for $f(z)=z^m$, $m\in\N$, $f(z)=\log{}z$,
has been proven in \cite{GO3}. 
%The computation of the
%coefficients is however a different problem. 
Explicit expressions for the first two coefficients
were given in \cite{GO3,GO}
in the case when $f(z)=z^m$, $m\in\N$, $f(z)=\log{}z$,
and for general $f\in\an$,
%by A.~Laptev, D.~Robert and Yu.~Safarov
in \cite{LRS}.
See Remark~\ref{CH1_s5rem1}
in subsection \ref{subsecRR} for a discussion
of the formulas for further coefficients. 
\end{rem}
\begin{rem}
\label{rem_RR}
The formula \eqref{CH1_s5ieFprin} for $W_3[f]$ 
has a structure similar to the
coefficient in the third asymptotic term in a Szeg\"o
type expansion for convolution operators obtained by
R.~Roccaforte in \cite{RR}. A certain combinatorial identity
known as {\em Spitzer's formula\ }(see \cite[Section~3.7]{DM}),
which is a version of the Bohnenblust--Spitzer 
combinatorial theorem, is also used in the proof
in \cite{RR}. (See also \cite{B1,Ro} and \cite[Remark~1.6]{GLc}).
%In Remark~\ref{CH1_s5rem1}
%in subsection \ref{subsecRR} we explain how
%the approach of \cite{RR} could be used in our setting. 
A second order Szeg\"o type expansion for convolution operators
was established by H.~Widom in \cite{Wn} %for convolution operators,
with the help of the usual Hunt--Dyson combinatorial
formula \eqref{CH1_s2e10orig}. A full asymptotic expansion
for convolution operators is obtained in \cite{Wbook}.
\end{rem}
%
%
%We explain a method of computing higher order terms in 
%Section \ref{CH1_s5} (see \eqref{CH1_s5e1var} and 
%Remark \ref{CH1_s5rem100}).
%
%The combinatorial formula gHD is presented in Section \ref{CH1_s4nn}. 
%
%
\begin{rem}
\label{CH1_s5iREMwz}
If $f\in\an$ is a polynomial then each of the sums
over $j,k,l$ in $\UpsilonPoi[f](B)$ reduces to a finite sum for
any Zoll manifold $M$ and any given $B\in\Psi^0(M)$, because
$\Phi_m[z^n]=0$ identically for $m>n$.
The term $\UpsilonPoi[f](B)$
vanishes for any $f\in\an$ in the case of $d=1$, 
and also in the following case:
Assume $f(z)=\sum_{k=1}^{\infty}c_{2k-1}z^{2k-1}$, 
let $M=\S^d$, $d\in\N$, $d\geq2$, and
let $B\in\Psi^0(M)$ have a principal symbol $b_0\in{}C^\infty(\S^d)$ which is 
independent of $\xi$, 
and in addition is an {\em odd\ } function on $\S^d$. 
Thus $b_0(x)=-b_0(x_a)$ for all pairs
of antipodal points $x,x_a\in\S^d$.
Then the Poisson brackets in $\UpsilonPoi[f](B)$
are odd with respect to $\xi$, and vanish after the 
integration over $S^*M$.
(We refer to \cite{W2} for an auxiliary calculation of
the needed Poisson bracket on $\S^2$, which can be easily 
modified for any $\S^d$, $d\geq3$.)
Under these assumptions the coefficient
of the first asymptotic term in \eqref{eq_th5} of order ${}n^d$, 
which is a part of $\Tr P_n{}f(B)P_n$, 
also vanishes, but
not the coefficients of the second and the third asymptotic terms
of orders $n^{d-1}$ and $n^{d-2}$, respectively. 
\end{rem}
\begin{rem}
In order to {\em define\ }the third Szeg\"o asymptotic term in \eqref{eq_th5}
the condition $f\in{}C^4(I)$ is necessary, where $I$ is the
closed set of values of $\sigma_0(B)$.
If the %non-symmetric
term $\UpsilonPoi[f](B)$
is absent (see Remark \ref{CH1_s5iREMwz}) then the condition $f\in{}C^3(I)$
is necessary. It was shown %by A.~Laptev, D.~Robert and Yu.~Safarov
in \cite{LRS} that if $B$ is self-adjoint and $f^{\prime\prime}\in{}L^\infty(I)$ 
then the second order Szeg\"o
formula holds, see also \cite{LS}. 
It would be interesting to know if the condition
of essential boundedness of the fourth derivative of $f$
%(or even $f\in{}C^4(I)$)
is sufficient for Szeg\"o asymptotic formula with three terms to hold.
Another question concerns the best possible norm
in place of \eqref{eq_seminorm}.
\end{rem}
\begin{rem}
\label{rem_comb}
We have discovered and proved
the generalized Hunt--Dyson combinatorial formula
(Theorem~\ref{CH1_S2TH1}) being unaware of the Bohnenblust--Spitzer
theorem (Theorem~\ref{CH1_s1COR2}). A derivation of the gHD from the BSt
and vice versa can be found in \cite{GLc}, see \cite[Chapter~2]{GiPhD}
for an independent proof of the gHD.
Also the importance of the BSt in the theory of the maximum
of a random walk with real i.i.d. steps and related results 
are discussed in \cite[Remark~1.5 and 1.6]{GLc}.
\end{rem}
\subsection{Explicit asymptotic formulas 
for $\log\det{}P_nBP_n$, as $n\ra\infty$}
%We would like to present explicitly the coefficients
%in the asymptotic expansion of $\log\det{}P_nBP_n$,
%as $n\ra\infty$, for $B\in{}\Psi^0(M)$.
Let $f(z)=\log z$.
Theorem~\ref{thm_onedim} and \ref{thm_main} give an expression
for $\log\det{}P_nBP_n=\tr\,\log{P_nBP_n}$ 
as a sum of $\tr{}P_n(\log{B})P_n$
and two lower order corrections, as $n\ra\infty$.
Proposition~\ref{prop3} below gives, for an arbitrary $G\in{}\Psi^0(M)$,
an auxiliary 
asymptotic expansion for 
%\begin{equation}
%%\label{eqstrace}
$\tr\,P_nGP_n$, as $n\ra\infty$.
%\end{equation}
%
%To prove Proposition~\ref{prop3}, 
%we sum over $k=1,\cdots,n$ in \eqref{eq_CdV}, see \cite{G1} for details.
For dimension $d=1,2$, we need the constant coefficient
in this expansion, which is more complicated than the other ones,
see the proof of Proposition~\ref{prop3} in Section~\ref{secpfprop3}
for details. 
Let $R_l(G)$, $l=0,1,2,\cdots$, be the Guillemin--Wodzicki residues
as given in \eqref{eqR} below.
%The terms of all orders
%in \eqref{eq_CdV}, and also the possible rapidly decaying term, 
%will contribute to it. 
%Let us therefore for $d=1,2$ 
For $d=1,2$ make an additional assumption
\begin{equation}
\label{eq_assu}
   \sum_{l=0}^\infty{}\big|R_l(G)\big|<\infty
\end{equation}
under which 
the following sums are absolutely convergent
\begin{equation}
\label{eq_C}
   C(G):= \sum_{k=1}^\infty \Big(\tr(\pi_kG)
   -\sum_{l=0}^{+\infty} k^{d-1-l}\,R_l(G)\Big),
\end{equation}
see Section~\ref{secpfprop3} for the proof.
%Note that in \eqref{eq_C} and also in the 
%proposition below there only appears 
%the series $\sum_{l=0}^\infty{}R_l(G)$.
%However we need the absolute convergence \eqref{eq_assu}
%in the proof of the remainder estimate.
 
Let $\gamma$ denote the Euler constant and $\zeta$ the Riemann
zeta function.
\begin{prop}
\label{prop3}
Let $M$ be a Zoll manifold of dimension $d\in\N$.
Let $P_n$ be as above and
assume that $G\in\Psi^0(M)$. For $d=1,2$, assume in addition
that \eqref{eq_assu} holds, and let $C(G)$ be defined by \eqref{eq_C}. 
Then the following holds as $n\ra\infty$,\newline
(i) for $d=1$,
$$
\begin{aligned}
 \tr\,&P_nGP_n = n\cdot R_0(G) + \log n\cdot R_1(G)\\
 &+ 
  \Big(C(G)+\gamma\,R_1(G)+\sum_{l=2}^\infty\zeta(l)\,R_{l}(G)\Big)
+\frac1n\cdot\Big(\frac12\,R_1(G) - R_2(G)\Big)  + O\Big(\frac1{n^{2}}\Big),
\end{aligned}
$$
(ii) for $d=2$,
$$
\begin{aligned}
 \tr\,P_nGP_n = n^2\cdot\frac12\,R_0(G) &+ 
  n\cdot\Big(\frac12\,R_0(G) + R_1(G)\Big)
+ \log n\cdot R_2(G)\\
 &+ 
  \Big(C(G)+\gamma\,R_2(G)+\sum_{l=2}^\infty\zeta(l)\,R_{l+1}(G)\Big)
  + O\Big(\frac1n\Big),
\end{aligned}
$$
(iii) for $d\geq3$,
$$
\begin{aligned}
 \null&\tr\,P_nGP_n = n^d\cdot\frac1d\,R_0(G) + 
  n^{d-1}\cdot\Big(\frac12\,R_0(G) + \frac1{d-1}\,R_1(G)\Big)\\
&+ 
  n^{d-2}\cdot\Big(\frac{d-1}{12}\,R_0(G) + \frac12\,R_1(G)
            + \frac1{d-2}\,R_2(G)\Big)
+ \log n\cdot R_d(G)
  + O(n^{d-3}).
\end{aligned}
$$
\end{prop}
\begin{rem}
\label{remtriv}
The coefficients of $n^d$, $n^{d-1}$
and $\log{}n$ for $d\geq2$ can be found in
\cite[after Lemma~0.2]{GO}.
From Proposition~\ref{prop3} with $G=\log{}B$
we see that $\tr{}P_n(\log{}B)P_n$
in Theorem~\ref{thm_onedim} and \ref{thm_main}
contributes to the leading asymptotic term of order $n^d$, and also to
all lower order terms of order $n^{j}$, $j=d-1,\cdots,1,0,-1,\cdots$,
and to the logarithmic term $\log{n}$, as $n\ra\infty$.
In the classical SSLT the situation is much simpler:
$\log{}B$ is just the Toeplitz matrix of the operator of multiplication
by $\log{}b$, and so $\tr{}P_n(\log{}B)P_n=(2n+1)\widehat{(\log{}b)}_0$.
\end{rem}
Now we are ready to state two corollaries. 
\begin{cor}
\label{cor4}
Let $B\in\Psi^0(\S^1)$ have a strictly positive principal
symbol and with a sufficiently small symbolic norm
$|||I-B|||_1$.
Assume also that \eqref{eq_assu} holds.
Then the following holds as $n\ra\infty$,
$$
 \log\det{}P_nBP_n = c_1\cdot n 
  +  c_{{\rm log}}\cdot\log n
          + c_0 + c_{-1} \cdot \frac1n + O\Big(\frac{1}{n^2}\Big),
$$
where the coefficients are the sums of the corresponding
coefficients from Theorem~\ref{thm_onedim} 
and Proposition~\ref{prop3}(i) for $G=\log{}B$.

Assume further that 
$\sigma_0(B)$ 
and $\sub(B)$ do not depend on the direction of $\xi$,
that is $\sigma_0(B)(x,\xi)=b_0(x)$ 
and $\sub(B)(x,\xi)=b_{{\sub}}(x)\,|\xi|^{-1}$, for $(x,\xi)\in{}S^*\S^1$.
Assume also that $b_{-2}=0$. Then the following holds as $n\ra\infty$,
\begin{equation}
\label{eqmain}
\begin{aligned}
      \log&\det{}P_nBP_n = n\cdot
     2 \int_0^{2\pi}\log b_0(x)\,\frac{dx}{2\pi}\\
                       &+ \log n\cdot2\int_0^{2\pi}
       \frac{b_{{\rm sub}}(x)}{b_0(x)}\,\frac{dx}{2\pi}\\ 
                      &+\Bigg(
    \sum_{k=1}^\infty k\,\widehat{(\log b_0)}_k
   \widehat{(\log b_0)}_{-k} 
+ C(\log{}B)+ \gamma\,R_1(\log{}B)+\sum_{l=2}^\infty\zeta(l)\,R_{l}(\log{}B)\Bigg)\\
         &+ \frac1n\cdot \Bigg(    \sum_{k=1}^\infty k\,
  \widehat{(\log b_0)}_k
\big(\widehat{b_{{\rm sub}}/b_0}\big)_{-k}
+\int_0^{2\pi}\bigg[\frac{b_{{\rm sub}}(x)}{b_0(x)}
                           +\bigg(\frac{b_{{\rm sub}}(x)}{b_0(x)}\bigg)^2\,\bigg]
         \,\frac{dx}{2\pi}\Bigg)\\
            &+O\bigg(\frac1{n^2}\bigg),
\end{aligned}
\end{equation}
where $C(\log{}B)$ is given by \eqref{eq_C}.
\end{cor}
The proof of \eqref{eqmain}
is an exercise in the calculus of PsDO's together with Proposition~\ref{prop3}(i),
and is left to the reader.
%\begin{rem}
In some simple cases, for instance for $b_{{\rm sub}}(x)=\pm\frac12{b_0(x)}$,
the left-hand side in \eqref{eqmain} can be computed explicitly.
The coefficients of $n$, $\log{}n$, and $\frac1n$ on the right in \eqref{eqmain}
in these cases are as expected,
%in the two simple cases $b_{{\rm sub}}(x)=\pm\frac12{b_0(x)}$,
see also Remark~\ref{remmatrix} below.
%\end{rem}
%
\begin{cor}
\label{cor5}
Let $M$ be a Zoll manifold of dimension $d\geq2$.
Assume that $P_n$ and $A$ are as in Theorem~\ref{thm_main}.
Let $B\in\Psi^0(M)$ have a strictly positive principal
symbol and a sufficiently small symbolic norm $|||I-B|||_1$.
For $d=2$, assume in addition \eqref{eq_assu}.
Then the following holds, as $n\ra\infty$,
\begin{equation}
\label{eq_expl}
 \log\det{}P_nBP_n = C_d^{(d)}\cdot n^d  + C_{d-1}^{(d)}\cdot n^{d-1}  
   +C_{d-2}^{(d)}\cdot n^{d-2}  
  +  C_{{\rm log}}^{(d)}\cdot\log n
          + O\big(n^{d-3}\big),
\end{equation}
where the coefficients are the sums of the corresponding
coefficients from Theorem~\ref{thm_main} for $f(z)=\log{}z$
and Proposition~\ref{prop3}(ii) or (iii) for $G=\log{}B$.
If one counts the logarithmic term, 
this expansion is fourth order for $d=2,3$
and third order for $d\geq4$. 
\end{cor}
\begin{rem}
\label{remdet}
%The {\em existence\ }of a full expansion of the type \eqref{eq_expl}
%has been proven in \cite{GO3}. The computation of the
%coefficients is a different problem.
The coefficients $C_d^{(d)}$ and $C_{d-1}^{(d)}$, $d\in\N$,
have been found in \cite{GO3,GO}.
The most interesting coefficient in \eqref{eq_expl}
is the constant one, since one can think
of $\exp{}C_0^{(d)}$ as a regularized determinant of $B\in\Psi^0(M)$,
see \cite[after~(3)]{GO3} and \cite{O2}. The sum 
$$
   \gamma{}R_d(\log{}B)+\sum_{l=2}^\infty\zeta(l)R_{l+d-1}(\log{}B)
$$
will for all $d\in\N$ be a part of $C_0^{(d)}$.
For $d=1$, Corollary~\ref{cor4} gives a full expression for $C_0^{(1)}$.
For $d=2$, Corollary~\ref{cor5} gives a full expression for $C_0^{(2)}$,
which is quite lengthy.
\end{rem}
\begin{rem}
Let us compare the result of Corollary~\ref{cor4} with
a generalization of SSLT to the case of $B$ being an operator
of multiplication by a function $b(x)$ having discontinuities
which is due to H.~Widom and E.~Basor.
In this case $\log{}b(x)$ also has discontinuities, and so the
series $\sum_{k\in\Z}|k|\,|\widehat{(\log{b})}_k|^2$
diverges logarithmically.
The following
third order asymptotic formula holds for 
the operator of multiplication by a piecewise $C^2$
function $b(x)$
\begin{equation}
\label{eqBasor}
       \log\det P_nBP_n = a_1\cdot n+a_2\cdot\log n+a_3+o(1),\qquad n\ra\infty,
\end{equation}
where $a_1$ as in \eqref{eqmain},
the coefficient $a_2$ has been
computed by H.~Widom 
in \cite{W4}, and 
the constant term $a_3$ has been found by E.~Basor in \cite{B}. 
Note that the matrix $B$ in \eqref{eqBasor} is still Toeplitz,
the logarithmic order of the subleading term being
due to a slower decay of the Fourier coefficients of $b(x)$. 
In our case the matrix of the operator
$B\in\Psi^0(\S^1)$ is {\em not\ }Toeplitz (see Remark~\ref{remmatrix} below),
and the $\log{}n$ term comes from
the contribution of $\sub(B)$.

%We would like to mention that t
It would be interesting to find a compact formula
for the constant term in \eqref{eqmain}.
We mention that the constant $a_3$ in \eqref{eqBasor}
found in \cite{B} has a form similar to the one in \eqref{eqmain}.
It contains a ``finite'' term and and an infinite series of
certain integrals multiplied by the values of the Riemann zeta function at 
the points $3,5,\cdots$. Interestingly, an ``invariant'' form
of that series has been found in \cite{Wotaa}. It is written as a single
integral involving the function 
$$
       \Psi(x):=\frac{d}{dx}\,\log \Gamma(x).
$$
This gives the hope that a similar formula can be
found for the constant in \eqref{eqmain}.
%in the infinite series in the constant in
%\eqref{eqmain} can also be rewritten in a compact form.
%
%%We have done some 
%%computations trying to find the constant term in \eqref{eqmain},
%%and the function $\Psi(x)$ has been appearing there.
\end{rem}
\begin{rem} 
\label{remmatrix}
The matrix interpretation of Corollary~\ref{cor4} is as follows.
Assume for simplicity that $B\in\Psi^0(\S^1)$ 
is as in the second part of Corollary~\ref{cor4},
that is $\sigma_0(B)(x,\xi)=b_0(x)$ 
and $\sub(B)(x,\xi)=b_{{\sub}}(x)\,|\xi|^{-1}$, 
for all $(x,\xi)\in{}S^*\S^1$.
Assume also that $b_{-2}=b_{-3}=\cdots=0$.
Let $B_0$ and $B_{{\rm sub}}$ be the operators of
multiplication by $b_0$ and $b_{{\rm sub}}$, respectively.
Let $D$ be the linear operator in $L^2(\S^1)$ such that
$$
        De^{ikx}=\begin{cases}\frac1{|k|}e^{ikx},&|k|\geq1\cr0,&k=0.\end{cases}
$$
Note that this is not a differential, but rather a smoothing operator
of order $-1$.
There is known a correspondence between the classical PsDO's on the
circle and their discrete counterparts, see \cite{M} for details.
By that correspondence, the zeroth order PsDO $B$ we started with
equals $B_0+B_{{\rm sub}}D$. 
Introduce two Toeplitz matrices,
$\widehat{B}_0:=\{\widehat{(b_0)}_{j-k}\}_{j,k\in\Z}$
and $\widehat{B}_{{\rm sub}}:=
\{\widehat{(b_{{\rm sub}})}_{j-k}\}_{j,k\in\Z}$.
Set also $\widehat{D}:=\diag(\cdots,\frac{1}{3},\frac12,1,0,1,\frac12,\frac{1}{3},\cdots)$.
Then the matrix representation of $B_0+B_{{\rm sub}}D$
is $\widehat{B}_0+\widehat{B}_{{\rm sub}}\cdot{}\widehat{D}$.
Finally, set $\widehat{P}_n=\diag(\cdots,0,1,\cdots,1,0,\cdots)$ ($2n+1$ ones).
Then Corollary~\ref{cor4} gives a fourth order asymptotics for
the determinant of the truncated matrix 
$\widehat{P}_n\cdot(\widehat{B}_0+\widehat{B}_{{\rm sub}}\cdot{}\widehat{D})\cdot\widehat{P}_n$.

Now we can reformulate the question 
of finding the constant term in \eqref{eqmain} in purely matrix terms.
Drop the hats and the dots %in $\widehat{D}$, $\widehat{P}_n$ 
for brevity.
Let $C_1$ be a Toeplitz matrix that corresponds to
the operator of multiplication by $b_{{\rm sub}}/b_0$,
and let the matrix $D$ be as above.
Clearly, the matrices $C_1$ and $D$ do not commute.
Assume that the matrix $\log(I-C_1D)$
is well-defined.
The question is to compute the constant coefficient 
in $\tr P_n\log(I-C_1D)P_n$,
or which is the same, the constant coefficient in
$$%\begin{equation}
%\label{eqcc}
        \tr{}P_n\log(I-D^{1/2}C_1D^{1/2})P_n, \qquad n\ra\infty.   
$$%\end{equation}
As we have noticed in Remark~\ref{remtriv},
this question is trivial for a Toeplitz matrix $T$ 
in place of $D^{1/2}C_1D^{1/2}$.
\end{rem}
\subsection{A related result for the maximum of a random walk}
Let us explain how we can use the $j$-maps $\Phi_j$, $j\in\N$,
to write the bivariate characteristic function
of the maximum of a random walk and its position at a smaller time.
There is a lack of symmetricity in this problem, and this
case is not considered in \cite{S1}. Let $X_1,X_2,\cdots$ be independent 
real valued random variables which assume real values 
and have the same distribution density $\phi(X)$.
Assume for simplicity that $\phi$ is Schwartz class.
(The result below holds for much more general $\phi$,
e.g., for discrete random variables, if understood 
in the sense of distributions.)
%[Assume that all the moments of $X$ exist and are finite.]
In our case the characteristic function
$$
     \E\{e^{itX}\}:= \int_{-\infty}^\infty e^{itX}\,\phi(X)\,dX
$$
is well-defined and invertible as a Fourier transform.
Let $S_p:=X_1+\cdots+X_p$ be the position at time $p\in\N$
of the random walk
starting off at the origin. For any $q\in\N$ 
introduce the random variable
$$
%\begin{aligned}
     T_{p+q}%&:=M_{p+q}(X_1,\cdots,X_p,X_{p+1},\cdots,X_{p+q})\\
         := \max\big(0,S_1,\cdots,S_p,S_{p+1},\cdots,S_{p+q}\big),
%\end{aligned}
$$
which is the length of the maximal excursion to the right of $0$
during the time interval $0,1,\cdots,p+q$. Note that the time $p+q$
is strictly larger than $p$.
%
%We compute the bivariate characteristic function for the random
%variables $S_p$ and $T_{p+q}$ for all
%$p,q\in\N$. We express the answer in form of a formal power series
%(generating function). Also we calculate explicitly the
%moment of $T_m$, $m\in\N$, of any order $n\in\N$.
%
Let us introduce a non-negative valued function
similar to the functions \eqref{CH1_s5ieOj1}, \eqref{CH1_s5ieOj2}
and \eqref{CH1_s5ieOj3}.
For arbitrary $j,k\in\N$ and $y_1,\cdots,y_j,z_1,\cdots,z_k\in\R$ we set
\begin{equation}
\label{CH1_s5iePOj2}
\begin{aligned}
      \Omega_{j,k}(&y_1,\cdots,y_j,z_1,\cdots,z_k)
          := (y_1)_++\cdots+(y_j)_+ \\
          &+\big(\,-(y_1)_--\cdots-(y_j)_-
           +(z_1)_++\cdots+(z_k)_+\,\big)_+.
\end{aligned}
\end{equation}
%Let the functors $\Phi_j$, $j\in\N$, be determined 
%by the conditions \eqref{CH1_s5iePhilin}
%and \eqref{CH1_s5iePhi2}.
Denote
$$
  \hat\phi(\eta):= \int_{-\infty}^{\infty} e^{-i\eta{}y}\,\phi(y)\,dy,\qquad
%$$
%and
%$$
   \F_{\eta\ra{}y}^{-1}[f(\eta)] := (2\pi)^{-1}
          \int_{-\infty}^{\infty} e^{iy\eta}\,f(\eta)\,d\eta.
$$
%We are ready to formulate the result for random walks.
\begin{thm}
\label{CH1_s5ith1RW}
%Assume $p,q\in\N$.
Let $\phi$, $X_p$, $S_p$, and $T_{p+q}$, $p,q\in\N$, be as above.
Let $\Phi_j$, $j\in\N$, be as in \eqref{CH1_s5iePhi2}.
Then for $|a|,|b|<1$ the following holds
\begin{equation}
\label{CH1_s5ie1RW}
\begin{aligned}
      \null\sum_{p,q=1}^\infty &a^pb^q\,
       \E\Big\{e^{i\alpha{}S_p+i\beta{}T_{p+q}}%m_{p+q}^*(X_1,\cdots,X_{p+q})}
\Big\}\\
     &=\sum_{j,k=1}^\infty \frac1{j!}\frac1{k!}
       \int_{-\infty}^\infty\!\cdots\! 
        \int_{-\infty}^\infty e^{i\alpha{}(y_1+\cdots+y_j)
       +i\beta\Omega_{j,k}(y_1,\cdots,y_j,z_1,\cdots,z_k)}\\
&\quad\quad\times\F_{ \genfrac{}{}{0pt}{}{\eta_1\ra{}y_1,\cdots,\eta_j\ra{}y_j}
           {\zeta_1\ra{}z_1,\cdots,\zeta_k\ra{}z_k} }^{-1}
\bigg[\Phi_j\Big[\frac{as}{1-as}\Big]
     \big(\hat\phi(\eta_1),\cdots,\hat\phi(\eta_j)\big)\\
&\quad\quad\quad\quad\times\Phi_k\Big[\frac{bt}{1-bt}\Big]
     \big(\hat\phi(\zeta_1),\cdots,\hat\phi(\zeta_k)\big)\bigg]
      \,dy_1\cdots dy_j\,dz_1\cdots dz_k.
\end{aligned}
\end{equation}
In this formula, the 
$j$-maps $\Phi_{j}$ and $\Phi_k$ act on the function of $s$ and $t$,
respectively. The arguments of these actions are the 
values of $\hat\phi$ at the corresponding points.
\end{thm}
The proof for fixed $p,q\in\N$ is carried out analogously to
the computations in Section~\ref{CH1_s6Poi}, 
see \cite[Section~1.10]{GiPhD} for details.
%The reason that 

%
\subsection{Organization of the article}
The paper is organized as follows.
In Section~\ref{CH1_s95}, we recall the method of \cite{GO},
estimate the remainder after the third Szeg\"o
term, and justify the passage
from the set of polynomials to an analytic function $f$.
After that we deal with an arbitrary 
monomial $f(z)=z^m$, $m\in\N$.
An operator $B_{\k_1}\cdots{}B_{\k_m}\in\Psi^0(M)$, 
$\k_1,\cdots,\k_m\in\Z$,
arises in that computation (see \eqref{eq_bk} below for the definition
of the Fourier coefficient $B_\k$, $\k\in\Z$).
In Section~\ref{CH1_s6aux}, we compute the contribution
of $\sigma_0\big(B_{\k_1}\cdots{}B_{\k_m}\big)$,
$m\geq2$, to higher order Szeg\"o terms.
The resulting expression involves only $\sigma_0(B)$.
In Section~\ref{CH1_s6sub}, we calculate the contribution of the
symmetric part of $\sub\!\big(B_{\k_1}\cdots{}B_{\k_m}\big)$,
$m\geq2$, to the third Szeg\"o term.
It depends on $\sigma_0(B)$ and $\sub(B)$.
In Section~\ref{CH1_s6Poi}, we compute the contribution of
the non-symmetric part of $\sub\!\big(B_{\k_1}\cdots{}B_{\k_m}\big)$,
$m\geq2$,
to the third Szeg\"o term, which involves Poisson brackets of the 
principal symbol of $B$ shifted along by geodesic flow. 
This contribution
depends only on $\sigma_0(B)$.
%%%%%, but is lengthy. It vanishes in certain special cases.
%
In Section~\ref{secpfprop3}, a proof of Proposition~\ref{prop3} is given.
In Section~\ref{CH1_s6Phi}, we find an expression 
for $\Phi_j[f]$, $j\in\N$, and also for $W_3[f]$, for an arbitrary $f\in\an$.
In Section~\ref{CH1_s4nn}, the auxiliary combinatorial background %material
is given. 

The results of this paper (Theorem \ref{thm_onedim}, \ref{thm_main} 
and \ref{CH1_S2TH1}, Proposition~\ref{prop3}
and Corollary~\ref{cor4} and \ref{cor5})
were announced in \cite{GL}, where we also gave an outline of the proofs.

{\bf Acknowledgments.}
This paper is a modified part of the author's 
Ph.D. thesis~\cite{GiPhD}.
I would like to express
a deep gratitude to my thesis adviser Ari Laptev, for suggesting
the problem and his constant attention to the work.

I would like to thank Percy Deift 
for giving a reference \cite{RS}
from which I learned about \cite{S1} and
the Bohnenblust--Spitzer
combinatorial theorem, shortly after having discovered
and proved the generalized Hunt--Dyson formula. 

I would also like to thank Kurt Johansson for suggesting
a derivation of the formula for 
$\tilde\Phi_j[f]$, $f\in\an$, in Section~\ref{CH1_s6Phi},
via generating functions.

Finally, I would like to thank the Department of Mathematics
of the Royal Institute of Technology (KTH), Stockholm,
for %excellent working conditions and 
a generous financial support
during the whole period of my graduate studies. 
\section%[The full third Szeg\"o term and the remainder estimate]
{The full third Szeg\"o term and the remainder estimate}
%{Proof of Theorem \ref{CH1_s5ith1Sz}: 
%Calculation of the full third Szeg\"o term
%and the remainder estimate}
\label{CH1_s95}
\subsection{The abstract scheme from \cite{GO}}
We start by expanding the analytic function $f(z)$
in a power series about $0$
and proving \eqref{eq_th5} for $f(z)=z^m$
for an arbitrary $m\in\N$. After that we justify the
passage from the set of polynomials to analytic functions.
Let us recall the method of \cite{GO}. 
Let $\pi_k$, $k\in\N$, be the projection on the $k$th eigenspace
of the operator $A$ and set $\pi_k:=0$ for $k\leq0$.
Then $P_n=\sum_{k=1}^n\pi_k$ for $n\in\N$, and we set
$P_n:=0$, $n\leq0$.
For an arbitrary $B\in\Psi^0(M)$ and $t\in\R$ 
introduce the operator
$$%\begin{equation}
%\label{eq_Bt}
       B^t:=e^{-itA}Be^{itA}.
$$%\end{equation}
By Egorov's theorem, $B^t\in\Psi^0(M)$, and also 
\begin{equation}
\label{eqEgo}
     \sigma_0(B^t)(x,\xi)=\sigma_0(B)(\Theta^t(x,\xi)),
\end{equation}
where $\Theta^t$ stands for the shift by $t$ units 
along the geodesic flow.
Note that because $\sub(A)=\const$ 
the following also holds \cite{Gui}
\begin{equation}
\label{eqEsub}
     \sub(B^t)(x,\xi)=\sub(B)(\Theta^t(x,\xi)).
\end{equation}
Because $\spec(A)=\N$, the operator $B^t$ is periodic
in $t$ with period $2\pi$.
Therefore we can introduce the Fourier expansion
$
          B = \sum_{{\k}\in\Z} B_{{\k}}
$
where $B_{{\k}}\in\Psi^0(M)$, ${\k}\in\Z$, is defined by
\begin{equation}
\label{eq_bk}
           B_{\k} = \sum_{k=1}^\infty 
                                     \pi_{k+{\k}}\,B\,\pi_k
                 =  \frac1{2\pi}\int_0^{2\pi} e^{i{\k}{}t}\,
                                  e^{-itA}Be^{itA}\,dt.
\end{equation}
This together with \eqref{eqEgo} and \eqref{eqEsub} implies
\begin{lem}
\label{lem0sub}
For any $B\in\Psi^0(M)$ and $\k\in\Z$
$$%\begin{equation}
%\label{eqEgok}
     \sigma_0(B_\k)(x,\xi)=\int_0^{2\pi} e^{i{\k}{}t}\,
                                  \sigma_0(B)(\Theta^t(x,\xi))\,\frac{dt}{2\pi}
$$%\end{equation}
and
$$%\begin{equation}
%\label{eqEsubk}
     \sub(B_\k)(x,\xi)=\int_0^{2\pi} e^{i{\k}{}t}\,
                                  \sub(B)(\Theta^t(x,\xi))\,\frac{dt}{2\pi}.
$$%\end{equation}
\end{lem}
For $m\in\N$ and ${\k}_1,\cdots,{\k}_m\in\Z$
introduce the notation
\begin{equation}
\label{eq_Mm}
  M_m(\overline{{\k}}):= \min(0,{\k}_1,{\k}_1+{\k}_2,\cdots,{\k}_1+\cdots+{\k}_m).
\end{equation}
Now using the remarkable commutation relation
\begin{equation}
\label{eq_cr}
      B_{\k} P_n = P_{n+{\k}} B_{\k},\qquad n\in\N,\,\, {\k}\in\Z,
\end{equation}
we move all the projectors to the left in the expression
$$
   (P_nBP_n)^m = \sum_{{\k}_1,\cdots,{\k}_m}P_nB_{{\k}_1}P_n
             B_{{\k}_2}P_n\cdots P_nB_{{\k}_m}P_n
$$ 
obtaining $P_n B^m P_n$ plus another term. This implies for all $n\in\N$
\begin{equation}
\label{CH1_s95e1}
\begin{aligned}
    \tr\,(&P_n B P_n)^m - \tr P_n B^m P_n\\
       &= -\sum_{{\k}_1+\cdots+{\k}_m=0}\,\tr\,\big(
          (P_n-P_nP_{n+{\k}_1}\cdots{}P_{n+{\k}_1+\cdots+{\k}_m}\big)
                B_{{\k}_1}\cdots B_{{\k}_m}\big)\\
&=-\sum_{{\k}_1+\cdots+{\k}_m=0}\,\,
          \sum_{j=M_m({\overline{{\k}}})+1}^0
               \tr(\pi_{n+k}\,B_{{\k}_1}\cdots B_{{\k}_m}).
\end{aligned}
\end{equation}
\begin{rem}The relation \eqref{eq_cr} follows
readily after writing both sides using the definition \eqref{eq_bk}.
We note that the Fourier coefficients $B_\k$
are PsDO's even if $B$ is a multiplication operator. Also
even in the simplest case $M=\S^1$, $B=b(x)$,
$
     B_\k\neq \hat{b}_\k e^{i\k x}.
$
This explains the %seeming 
non-symmetricity of \eqref{eq_cr}
with respect to positive and negative $\k$.
\end{rem}
Next, for any $G\in\Psi^0(M)$, $M$ being a Zoll manifold,
there exists a full asymptotic expansion for $\tr(\pi_k{}G)$, as $k\ra\infty$,
see Lemma~\ref{lem_CdV} below.
This result is due to Y.~Colin de Verdi\`ere \cite{CdV}.
The coefficients in that expansion
are certain Guillemin--Wodzicki residues.
Recall that for any compact closed manifold $M$
of dimension $d\in\N$
the Guillemin--Wodzicki residue 
of a pseudodifferential operator $G\in\Psi^m(M)$ of order $m\in\Z$
is defined by 
\begin{equation}
\label{CH1_s5e2}
   \Res(G):= \int_{S^*M} \sigma_{-d}(G)(x,\xi)\,dxd\xi
\end{equation}
(recall that we have included $(2\pi)^{-d}$ 
in the notation $dxd\xi$).
For an arbitrary $G\in\Psi^0(M)$ denote
\begin{equation}
\label{eqR}
%\begin{aligned}
      R_l(G) :=\Res (A^{-d+l}G),\qquad l=0,1,2,\cdots.
%\end{aligned}
\end{equation}
\begin{lem}
\label{lem_CdV}
Let $M$ be a Zoll manifold of dimension $d\in\N$.
Assume $G\in\Psi^0(M)$.
Then for any $N=0,1,2,\cdots$, 
there exists $C_N<\infty$ such that
\begin{equation}
\label{eq_CdV}
         \Big| \tr(\pi_kG)
   -\sum_{l=0}^{N} k^{d-1-l}\,R_l(G) \Big|\leq C_N\,k^{d-2-N}, \quad k\in\N.
\end{equation}
\end{lem}
%\noindent
See \cite{CdV} and \cite[Appendix]{GO} for the proof of \eqref{eq_CdV}.
We need the expressions for the first two residues.
\begin{lem} 
\label{lem_res}
Let $A$ be defined by \eqref{eq_A}.
For any $G\in\Psi^0(M)$
\begin{equation}
\label{CH1_s9e1.1}
%\begin{aligned}
      R_0(G) =  \int_{S^*M}
                    \sigma_0(G)\,dxd\xi
%\end{aligned}
\end{equation}
and
\begin{equation}
\label{CH1_s9e1.1sub}
\begin{aligned}
      R_1(G) =\int_{S^*M} \big(\,(d-1)(\alpha/4)\,\sigma_0(G) + \sub(G)\,\big)
                   \,dxd\xi.
\end{aligned}
\end{equation}
\end{lem}
\begin{proof}
The equality \eqref{CH1_s9e1.1} 
follows easily from the computation rules for PsDO's
%we note that $A^{-d}G\in\Psi^{-d}(M)$
%from the fact that $\sigma_{-d}(A^{-d}G)$ is the principal symbol of $A^{-d}G$,
%the computation rule 
%for the principal symbol $\sigma_{-d}(A^{-d}G)=(\sigma_1(A))^{-d}\sigma_0(G)$,
and the fact that $\sigma_1(A)=1$ on $S^*M$.
Next, the equality 
\begin{equation}
\label{eqzzz}
       R_1(G)
             = \int_{S^*M}\sub(A^{-d+1}G)
                   \,dxd\xi.
\end{equation}
follows from the definition of the subprincipal symbol
(see, e.g., \cite{DG}), and the fact that the integral
over $S^*M$ of a Poisson bracket equals zero.
(We could also refer to Proposition~29.1.2 in \cite{Ho4},
differentiate with respect to the spectral parameter,
make the same remark concerning
the integral involving the Poisson bracket.)
Now \eqref{CH1_s9e1.1sub} follows from \eqref{eqzzz}
by the computation rules concerning the subprincipal symbol:
%In our case, $A^{-d+1}\in\Psi^{-d+1}(M)$ and $G\in\Psi^0(M)$. 
by \cite[(1.4)]{DG} 
$$%\begin{equation}
%%\label{CH1_s9e3}
\begin{aligned}
  \sub(A^{-d+1}G) &= \sub(A^{-d+1})\,\sigma_0(G) 
         + \sigma_{-d+1}(A^{-d+1})\,\sub(G)\\
                     &\quad+ (2i)^{-1}\{\sigma_{-d+1}(A^{-d+1}),
                          \sigma_0(G)\},
\end{aligned}
$$%\end{equation}
where $\{\cdot,\cdot\}$ denotes the Poisson bracket. 
Now by \cite[(1.3)]{DG}
$$
   \sub(A^{-d+1}) = (-d+1)\,\big(\sigma_1(A)\big)^{-d}\,\sub(A),
$$
and we use the definition \eqref{eq_A},
%\begin{equation}
%     \sub(A^{-d+1}) = (-d+1)\,\big(\sigma_1(A)\big)^{-d}\,\sub(A).
%\end{equation}
%For the operator $A\in\Psi^1(M)$ defined by \eqref{eq_A} one has
%$$
%            \sub(A)=\frac12\,\frac{\sub(-\Delta)}{\sigma_1(A)}-\frac{\alpha}{4}
%                                     =-\frac{\alpha}{4},
%$$
and that $\sub(-\Delta)=0$ to complete the proof.
\end{proof}
\subsection{Beginning of the proof of Theorem \ref{thm_main}}
We do not assume that the Fourier expansion
$B=\sum_{\k\in\Z}B_\k$ has only a finite number of terms.
Denote 
$$
   \vka:=(\k_1,\cdots,\k_m),\qquad B_{\vka}:=B_{\k_1}\cdots B_{\k_m}.
$$
Below $c(B)$ and $C(B)$ will denote various constants
which depend only on $B$ (and not on $m$).
We will indicate 
which seminorm of $B$ enters a certain $C(B)$ when necessary.
For any $M>0$ denote
\begin{equation}
\label{s95eqABR}
          B^{(M)}:=\sum_{|\k|\leq{}M}B_\k.
\end{equation}
Introduce for $m\in\N$, $M>0$ the set of indices
$$
     Q^m(M):=\{(\k_1,\cdots,\k_m)\in\Z^m\,:\,
                       |\k_l|\leq{}M,\,l=1,\cdots,m\}.
$$
We split the right-hand side of \eqref{CH1_s95e1} 
into the following three sums. The first sum is
\begin{equation}
\label{CH1_s95e115}
\begin{aligned}
    -\sum_{ \genfrac{}{}{0pt}{}{\vka\in{}Q^m({}n/(2m))}
                     {\k_1+\cdots+\k_m=0} }\,\,
          \sum_{j=M_m(\vka)+1}^0
             \Big[ &({}n+j)^{d-1}\,\res(A^{-d}B_{\vka})\\
                    &+ ({}n+j)^{d-2}\,\res(A^{1-d}B_{\vka}) \Big],
\end{aligned}
\end{equation}
the second
\begin{equation}
\label{CH1_s95e117}
\begin{aligned}
    -&\sum_{ \genfrac{}{}{0pt}{}{\vka\in{}Q^m({}n/(2m))}
                     {\k_1+\cdots+\k_m=0} }\,\,
          \sum_{j=M_m(\vka)+1}^0
               \Big[\tr(\pi_{{}n+j}\,B_{\vka})\\
                      &\quad-({}n+j)^{d-1}\,\res(A^{-d}B_{\vka}) 
                  -({}n+j)^{d-2}\,\res(A^{1-d}B_{\vka})\Big]
\end{aligned}
\end{equation}
and the third
\begin{equation}
\label{CH1_s95e119}
\begin{aligned}
    -\sum_{ \genfrac{}{}{0pt}{}{\vka\in{}\Z^m\setminus Q^m({}n/(2m))}
                     {\k_1+\cdots+\k_m=0} }\,\,
          \sum_{j=M_m(\vka)+1}^0
               \tr(\pi_{{}n+j}\,B_{\vka}).
\end{aligned}
\end{equation}
The contributions to the second and third Szeg\"o term 
come from \eqref{CH1_s95e115}, the expressions
\eqref{CH1_s95e117} and \eqref{CH1_s95e119} having a lower order, as $n\ra\infty$. 
In the rest of this section 
we estimate the remainder after the third Szeg\"o term.
It is very important that the symbol of the operator $B$
of any order is smooth. 
We will repeatedly refer to subsection \ref{s95subsecpflem}
where the most technical part of the computation
is carried out.
\subsection{Computation of the sum \eqref{CH1_s95e115}}
\label{subsecRR}
Let us single out the terms involving ${}n^{d-1}$ and ${}n^{d-2}$
in \eqref{CH1_s95e115}. 
Then the latter splits into the sum of
\begin{equation}
\label{CH1_S95E1151}
\begin{aligned}
    - \sum_{ \genfrac{}{}{0pt}{}{\vka\in{}Q^m({}n/(2m))}
                     {\k_1+\cdots+\k_m=0} }\,\,
          &\sum_{j=M_m(\vka)+1}^0
             \bigg[ \,{}n^{d-1}\cdot\res(A^{-d}B_{\vka})\\
                    &+ {}n^{d-2}\cdot\Big(\,(d-1)\cdot{}j\cdot\res(A^{-d}B_{\vka})
                         +\res(A^{1-d}B_{\vka}\,)\Big) \bigg]
\end{aligned}
\end{equation}
and the remainder 
\begin{equation}
\label{CH1_s95e1152}
\begin{aligned}
    - \sum_{ \genfrac{}{}{0pt}{}{\vka\in{}Q^m({}n/(2m))}
                     {\k_1+\cdots+\k_m=0} }\,\,
          &\sum_{j=M_m(\vka)+1}^0
             \bigg[ \,\sum_{k=0}^{d-3}{}n^k\Big(\genfrac{}{}{0pt}{}{d-1}{k}\Big)
j^{d-1-k}\cdot\res(A^{-d}B_{\vka})\\
                    &+ \sum_{k=0}^{d-3}{}n^k\Big(\genfrac{}{}{0pt}{}{d-2}{k}\Big)
j^{d-2-k}\cdot\res(A^{1-d}B_{\vka}) \bigg].
\end{aligned}
\end{equation}
%It will be shown in Section \ref{CH1_s9} that 
Taking into account Lemma \ref{lem_res} for $G=B_{\vka}$ and the formulas
$$
   \sum_{j=M_m({\overline{\k}})+1}^0 1 = -M_m({\overline{\k}}),\qquad
   \sum_{j=M_m({\overline{\k}})+1}^0 j  = - \frac{(\,M_m({\overline{\k}})\,)^2 
                           + M_m({\overline{\k}})}{2}
$$
we conclude that \eqref{CH1_S95E1151} equals
\begin{equation}
\label{CH1_s9e24}
\begin{aligned}
         \null &\sum_{ \genfrac{}{}{0pt}{}{\vka\in{}Q^m({}n/(2m))}
                     {\k_1+\cdots+\k_m=0} }  \bigg(  n^{d-1}\cdot{}
  M_m({\overline{\k}})\,\int_{S^*M}\sigma_0(B_{\vka})\,dxd\xi \\
&\qquad+n^{d-2}\cdot\bigg[\frac{d-1}{2}\Big(
   \big(M_m({\overline{\k}})\big)^2+\big(1+\frac{\alpha}{2}\big)\,M_m({\overline{\k}})\Big)\,
     \int_{S^*M}\sigma_0(B_{\vka})\,dxd\xi \\
&\qquad\qquad\qquad\qquad\quad\qquad
+M_m({\overline{\k}})\,\int_{S^*M}\sub(B_{\vka})\,dxd\xi
\bigg]\bigg).
 \end{aligned}
\end{equation}
In Section \ref{CH1_s6aux}, the expression
\begin{equation}
\label{CH1_s9.e5}
%\begin{aligned}
   \sum_{\k_1+\cdots+\k_m=0}\,
            (\,M_m({\overline{\k}})\,)^n \,
   \int_{S^*M}\sigma_0(G_{\k_1}\cdots{}G_{\k_m})\,dx\,d\xi
%\end{aligned}
\end{equation}
for $n\in\N$ and any $G\in\Psi^0(M)$ is computed (we set $G:=B^{({}n/(2m))}$).
The important point here is that both the domain of summation
and the second factor in \eqref{CH1_s9.e5} are symmetric,
and so we can symmetrize the first factor.
The integral in \eqref{CH1_s9.e5} is symmetric 
with respect to $\k_1,\cdots,\k_{m}$ because already the integrand is symmetric
$$
     \sigma_0(G_{\k_1}\cdots{}G_{\k_m}) = \sigma_0(G_{\k_1})\cdots{}\sigma_0(G_{\k_m}).
$$
The gHD for $n=1,2$ (Theorem~\ref{CH1_S2TH1})
is needed in this computation.

In Section \ref{CH1_s6sub} and \ref{CH1_s6Poi}, the expression
$$%\begin{equation}
%\label{CH1_s9.e5sub}
%\begin{aligned}
   \sum_{\k_1+\cdots+\k_m=0}\,
            M_m({\overline{\k}}) \,
   \int_{S^*M}\sub(G_{\k_1}\cdots{}G_{\k_m})\,dx\,d\xi
%\end{aligned}
$$%\end{equation}
for any $G\in\Psi^0(M)$ is computed (we set $G:=B^{({}n/(2m))}$).
Note that 
\begin{equation}
\label{CH1_s9e10}
\begin{aligned}
  \sub(G_{\k_1}\cdots{}G_{\k_m}) = &\sum_{k=1}^m\sub(G_{\k_k})
                       \prod_{ \genfrac{}{}{0pt}{}{p=1}{p\neq{}k} }^m \sigma_0(B_{\k_p}) \\
                       +\frac1{2i} &\,\sum_{1\leq{}k<l\leq{}m}
                               \{\sigma_0(G_{\k_k}),\sigma_0(G_{\k_l})\}\,
           \prod_{ \genfrac{}{}{0pt}{}{p=1}{p\neq{}k,p\neq{}l} }^m \sigma_0(G_{\k_p})
\end{aligned}
\end{equation}
where the first sum is symmetric with respect to $\k_1,\cdots,\k_{m}$,
and the second one is generally speaking not (not even after the integration
over $S^*M)$. Because of that circumstance the original method of symmetrization
\cite{K,S1,O,GO,GO2,O2} fails, and we have to modify it.
It turns out that each of the $m(m-1)/2$ terms
in the second sum in \eqref{CH1_s9e10} possesses a {\it partial\ }symmetry.
For instance, if 
$1\leq{}r$, $r+2\leq{}s$ and $s+2\leq{}m$ then in the expression
$$
%\begin{aligned}
       \sigma_0\big(G_{\k_1}\cdots G_{\k_r}\big)\,
      \big\{\sigma_0(G_{\k_{r+1}}),\sigma_0(G_{\k_{s+1}})\big\}\,
     \sigma_0\big(G_{\k_{r+2}}\cdots G_{\k_s}\big)\,
          \sigma_0\big(G_{\k_{s+2}}\cdots G_{\k_m}\big)\,
%\end{aligned}
$$
the indices $\k_1,\cdots,\k_r$ can be permuted without changing
the resulting expression, and the same holds for the groups
of indices $\k_{r+2},\cdots,\k_s$ and $\k_{s+2},\cdots,\k_m$
(and we are even allowed to interchange the indices within the three groups).
However we can neither interchange any index
from any of the three groups with $\k_{r+1}$, $\k_{s+1}$,
nor interchange the latter two. 
Here the original form of the BSt (Theorem~\ref{CH1_s1COR2}) is needed. 
\begin{rem}
\label{CH1_s5rem1} 
%
%In view of \eqref{CH1_s5e7} and \eqref{CH1_s5e7bis} and the latter fact,
In view of \eqref{CH1_s95e1} and Lemma \ref{lem_CdV},
 the computation of {\em all\ }terms in the Szeg\"o asymptotics
%formula \eqref{CH1_s5e4} 
is reduced
to the evaluation of the following expression
for different $n\in\N$ and $l=0,1,2,\cdots$,
\begin{equation}
\label{CH1_s5e1var}
       \sum_{\k_1+\cdots+\k_m=0}\,\,
          (\,M_m(\k_1,\cdots,\k_m)\,)^n\,
\res(A^{-d+l}B_{\k_1}\cdots B_{\k_m}).
\end{equation}
The gHD (Theorem \ref{CH1_S2TH1})
makes the computation possible for any $n\in\N$, 
provided that the second factor % Guillemin--Wodzicki residue 
in \eqref{CH1_s5e1var}
is symmetric with respect to $\k_1,\cdots,\k_m$.
The problem is that this is the case generally speaking only for $l=0$.
The case of $l=1$ and $n=1$ is dealt with in
Section \ref{CH1_s6sub} and \ref{CH1_s6Poi}.

We would like to 
mention that in the work \cite{RR} by R.~Roccaforte an interesting way to
rewrite expressions of the type \eqref{CH1_s5e1var}
with non-symmetric (matrix-valued) second factor has been suggested. 
The idea is to consider all the cases when the minimum is attained on
the sum
$$
       \k_1+\cdots+\k_j,\qquad j=1,\cdots,m-1,
$$
and make certain changes of summation indices. More precisely,
for an operator $B\in\Psi^0(M)$ written as a Fourier
series $B=\sum_{\k\in\Z}B_\k$ following \cite{RR} we set
$$%\begin{equation}
%\label{CH1_s5e300}
\begin{aligned}
       B_-&:=\sum_{\k<0}B_\k,\qquad B_-^{j+1}:= \big(B*B_-^{j}\big)_-,
              \quad j\in\N\\
       B_+&:=\sum_{\k\geq0}B_\k,\qquad B_+^{j+1}:= \big(B_+^{j}*B\big)_+,
    \quad j\in\N,
\end{aligned}
$$%\end{equation}
where $*$ denotes the discrete convolution on the Fourier series side. 
The order of the operators being convolved is important.
Then  \eqref{CH1_s5e1var} can be rewritten as 
\begin{equation}
\label{CH1_s5e1var2}
       \sum_{\k}\,\,
          (-(\k)_-)^n\,\res\bigg[A^{-d+l}\,
             \sum_{j=1}^{m-1}\Big(B_-^j\Big)_\k\, 
                         \Big(B_+^{m-j}\Big)_{-\k}\bigg].
\end{equation}
Note that the summation in \eqref{CH1_s5e1var2} is 
over a {\em single\ }variable $\k$.
This reminds the expression which appears
in the usual second Szeg\"o term \cite{K,GO}.
Also the formulas from \cite{W} have a similar structure.
However the formula \eqref{CH1_s5e1var2} for 
$f(z)=z^m$ is not very explicit,
and we could not write a reasonable formula for
an arbitrary analytic $f(z)$ starting from \eqref{CH1_s5e1var2}
(nor was it done in \cite{RR}).
\end{rem}
Bringing together the results of Section 3, 4, and 5,
we conclude that \eqref{CH1_s9e24} equals
\begin{equation}
\label{eq_together}
      {}n^{d-1}\cdot\Upsilondva[z^m]\big(B^{({}n/(2m))}\big) 
             + {}n^{d-2}\cdot\Upsilontri[z^m]\big(B^{({}n/(2m))}\big).
\end{equation}
This in turn equals
$$%\begin{equation}
%%\label{CH1_s95e1153}
 {}n^{d-1}\cdot\Upsilondva[z^m](B) 
  + {}n^{d-2}\cdot\Upsilontri[z^m](B)
$$ %\end{equation}
plus an error whose absolute value does not exceed
\begin{equation}
\label{CH1_s95e1153}
\begin{aligned}
       n^{d-1}\cdot&\Big|\,\Upsilondva
                              [z^m]\big(B^{({}n/(2m))}\big) 
               - \Upsilondva[z^m](B)\,\Big|\\
     + &n^{d-2}\cdot\Big|\,\Upsilontri
                              [z^m](B^{({}n/(2m))}) 
               - \Upsilontri[z^m](B)\,\Big|,\qquad n\in\N.
\end{aligned}
\end{equation}
We need the following auxiliary statement which is proved in
subsection \ref{s95subsecpflem}.
\begin{lem}
\label{CH1_s95l0} 
For $l=2,3$, an arbitrary $B\in\Psi^0(M)$ and any $L\in\N$ 
there exist two constants
$c_{l,L}(B)$ and $C_l(B)$ such that for $m\in\N$, $m\geq2$
$$%\begin{equation}
%\label{CH1_s95e0}
\begin{aligned}
    \Big|\,\Upsilon_l[z^m](B^{(M)}) 
               &- \Upsilon_l[z^m](B)\,\Big|\\
        &\leq c_{l,L}(B)\cdot m^2\cdot(C_l(B))^m\cdot M^{-L},\quad M\geq1.
\end{aligned}
$$%\end{equation}
\end{lem}
%Now estimate for instance the first term in \eqref{CH1_s95e1153}.
It follows from the proof of Lemma \ref{CH1_s95l0}
that $C_l(B)$, $l=2,3$, 
involves $\|\nabla^{6}\sigma_0(B)\|_\infty$
and $\|\nabla^2\sub(B)\|_\infty$, but does not depend on $L$.
 Applying Lemma \ref{CH1_s95l0} with $M={}n/(2m)$
and $l=2$, $L=2$, respectively, $l=3$, $L=1$ to the first,
respectively, second term in \eqref{CH1_s95e1153},
we estimate the latter by
$$
        {}n^{d-3}\cdot{}c\cdot m^4\cdot c(B)\cdot(C(B))^m,\quad n\in\N.   
$$
where $C(B)$, $l=2,3$, depends on $\|\nabla^{6}\sigma_0(B)\|_\infty$
and $\|\nabla^2\sub(B)\|_\infty$.

We finish this subsection by showing that
the absolute value of the
remainder \eqref{CH1_s95e1152} can be estimated, for $n\in\N$, by
$$
        {}n^{d-3}\cdot c(B)\cdot m\cdot(C(B))^m.
$$
For instance, the absolute value of
each of the $d-2$ terms in the first sum in
the square bracket in \eqref{CH1_s95e1152} does not exceed
\begin{equation}
\label{s95eqEXTRA}
\begin{aligned}
    {}n^{d-3}\cdot &c_d\,\sum_{ (\k_1,\cdots,\k_m)\in{}\Z^m }\,\,
              |M_m(\vka)|^{d}\,
                  \big|\,\Res(A^{-d}B_{\k_1}\cdots B_{\k_m})\,\big|\\
    &\leq{}n^{d-3}\cdot c_d\,\sum_{(\k_1,\cdots,\k_m)\in{}\Z^m}\,\,
              %\big(\,|\k_1|+\cdots+|\k_m|\,\big)^{d}\,\\
       (1+|\k_1|)^d\cdots(1+|\k_m|)^{d}\\
                  &\qquad\times\int_{S^*M} dx d\xi\, 
             \big|\sigma_0(B_{\k_1})\big|\cdots \big|\sigma_0(B_{\k_m})\big|\\
      &\leq{}n^{d-3}\cdot c(B)\cdot m\cdot(C(B))^m,\quad n\in\N.
\end{aligned}
\end{equation}
The last inequality is due to the
fast decay of the maximum over $S^*M$ of
$|\sigma_0(B_{\k})|$, as $|\k|\ra\infty$ 
(smoothness of $\sigma_0(B)$ and Egorov's theorem),
see subsection~\ref{s95subsecpflem} for details.
We remark also that in \eqref{s95eqEXTRA} $C(B)$ depends 
on $\|\nabla^{d+2}\sigma_0(B)\|_\infty$.
%\cite[Lemma~1.3]{GO} and Lemma~\ref{lem_rd} below.
In the same way one shows that the absolute value of
each of the $d-2$ terms in the second sum in
the square bracket in \eqref{CH1_s95e1152} is estimated by the
right-hand side of \eqref{s95eqEXTRA}.
\subsection{Computation of the sum \eqref{CH1_s95e119}}
For any $\k\in\Z^m\setminus{}Q^m({}n/(2m))$ the absolute value
of at least
one of the components of $\k$ is $>{}n/(2m)$. Redenote it by $\k_1$. 
We also note that $j\geq{}-n$ in \eqref{CH1_s95e119}
(otherwise $\pi_{{}n+j}=0$). Then the absolute value of 
the sum \eqref{CH1_s95e119} does not exceed
\begin{equation}
\label{CH1_s95e100}
\begin{aligned}
    m\sum_{j=-n}^0\,\, &\sum_{ \genfrac{}{}{0pt}{}{\k\in\Z^m:|\k_1|>{}n/(2m)}
                     {\k_1+\cdots+\k_m=0} }\,\,
               \Big|\tr(\pi_{{}n+j}\,B_{\k_1}\cdots B_{\k_m})\Big|\\
    &\leq m\sum_{j=-n}^0\,\, \sum_{ \genfrac{}{}{0pt}{}{\k\in\Z^m:|\k_1|>{}n/(2m)}
                     {\k_1+\cdots+\k_m=0} }\,\,
                c_d({}n+j)^{d-1}\,\|B_{\k_1}\|\cdots\|B_{\k_m}\|\\
    &\leq {}n^d\cdot m\cdot c_d \sum_{|\k_1|>{}n/(2m)}\|B_{\k_1}\|\,
                 \Big(\,\sum_{\k\in\Z}\|B_{\k}\|\,\Big)^{m-1},
\end{aligned}
\end{equation}
where the factor $c_d({}n+j)^{d-1}$ is an estimate 
of the multiplicity of the eigenvalue ${}n+j$
and $\|\cdot\|$ is the operator norm in $L^2(M)$. 
By \eqref{eq_norm0} below, there is a constant $c(B)$ such that
$$
        \sum_{|\k_1|>{}n/(2m)}\|B_{\k_1}\| 
        \leq c(B)\,\bigg(\frac{n}{2m}\bigg)^{-3},\quad n\in\N.
$$
Set $C(B):=\sum_{\k\in\Z}\|B_{\k}\|$.
By \eqref{eq_norm} below, $C(B)\leq\|B\|+(\pi^2/3)\big\|[A,[A,B]]\big\|$.
It follows now from \eqref{CH1_s95e100} that the absolute
value of the sum \eqref{CH1_s95e119} is estimated, for $n\in\N$, by
$$
       {}n^{d-3}\cdot c_d\,c(B)\cdot m^4\cdot (C(B))^m.
$$
\subsection{Computation of the sum \eqref{CH1_s95e117}}
For an arbitrary $G\in\Psi^0(M)$ and any $\mu\in\N$
denote
\begin{equation}
\label{s95eqTmu}
            T_\mu{}G:= \mu^{-d+3}\Big[\tr(\pi_{\mu}\,G)
                      -\mu^{d-1}\,\res(A^{-d}G) 
                  -\mu^{d-2}\,\res(A^{1-d}G)\Big].
\end{equation}
We wish to prove that there exist two constants $c(B)$ and $C(B)$
so that the limit, as $n\ra\infty$, of the absolute value of
\eqref{CH1_s95e117} divided by ${}n^{d-3}$ is $\leq{}c(B)\cdot m\cdot(C(B))^m$.
Let us estimate the absolute value of \eqref{CH1_s95e117} 
divided by ${}n^{d-3}$ by
\begin{equation}
\label{CH1_s95e500}
\begin{aligned}
    {}\sum_{ \genfrac{}{}{0pt}{}{\vka\in{}Q^m({}n/(2m))}
                     {\k_1+\cdots+\k_m=0} }\,\,
          &\sum_{j=M_m(\vka)+1}^0
               \Big(\frac{n+j}{n}\Big)^{d-3}\,\big|\,T_{{}n+j}B_{\vka}\,\big|\\
     &\leq\sum_{ \genfrac{}{}{0pt}{}{\vka\in{}Q^m({}n/(2m))}
                     {\k_1+\cdots+\k_m=0} }\,\,
          \sum_{j=M_m(\vka)+1}^0
               \,\big|\,T_{{}n+j}B_{\vka}\,\big|,
\end{aligned}
\end{equation}
where we have used that ${}n/2\leq{}n+j\leq{}n$ for $\vka\in{}Q^m({}n/(2m))$.
We need to interchange the $\lim_{{}n\ra\infty}$ to the right-hand
side in \eqref{CH1_s95e500} and the sums over infinitely growing
sets.
To justify that we find a summable over $\vka$ and $j$ majorant
and then refer to the Lebesgue dominated convergence.
We use the principle of uniform boundedness to present
a summable majorant. 

Introduce a Banach space $X=\Psi^0(M)/\Psi^{-3}(M)$ 
with a norm $\|\cdot\|_X$ given by
$$
\begin{aligned}
   \|G\|_X:= \|G\| + &\int_{S^*M}\big(\,
        |\sigma_{-d}(A^{-d}G)| 
      +|\sigma_{-d}(A^{1-d}G)| \\ 
      &+ |\sigma_{-d}(A^{2-d}G)|\,\big)\,dxd\xi, \quad G\in X,
\end{aligned}
$$
where the integrands are well-defined by the definition of
the Guillemin--Wodzicki residue \eqref{CH1_s5e2}.
Then it is clear from the definition \eqref{s95eqTmu}
that each $T_\mu:X\ra\C$ is a linear and bounded functional,
because the multiplicity of each eigenvalue $\mu$ is finite and 
$$
       |\res(A^{l-d}G)|\leq\int_{S^*M} \big|\,\sigma_{-d}(A^{l-d}G)\,\big|
                 \,dxd\xi,\quad l=0,1,2.
$$
Also by Lemma \ref{lem_CdV} for all $G\in{}X$ there exists the limit
\begin{equation}
\label{s95eqL31}
             \lim_{\mu\ra\infty}T_\mu G = \res(A^{2-d}G).
\end{equation}
Therefore for all $G\in{}X$ one has $\sup_{\mu\in\N}|T_\mu{}G|<\infty$.
By the principle of uniform boundedness, there exists a constant
$N<\infty$ such that
$$
           \sup_{\mu\in\N}|T_\mu{}G|\leq N\,\|G\|_X,\qquad G\in{}X.
$$
We finish the construction of the majorant by showing that it is summable.
The right-hand side in \eqref{CH1_s95e500} is estimated by 
$$
\begin{aligned}
      N\,\sum_{ \genfrac{}{}{0pt}{}{\vka\in{}Q^m({}n/(2m))}
                     {\k_1+\cdots+\k_m=0} }\,\,
          &\sum_{j=M_m(\vka)+1}^0
               \big\|\, B_{\vka} \,\big\|_X 
\leq   N\,\sum_{ \genfrac{}{}{0pt}{}{\vka\in{}Q^m({}n/(2m))}
                     {\k_1+\cdots+\k_m=0} }\,\,
          \big|M_m(\vka))\big|\,
               \big\|\, B_{\vka} \,\big\|_X \\
     &\leq   N\,\sum_{ \vka\in\Z^m }\,\,
          \big(\,|\k_1|+\cdots+|\k_m|\,\big)\,
               \big\|\, B_{\vka} \,\big\|_X
           < \infty,
\end{aligned}
$$
where the last inequality is due to 
the smoothness of the symbol of $B$ of any order,
%the rapid decay
%of the symbol of any order of $B_\k$, $\k\in\Z$, 
%Lemma~\ref{lem_rd}.
see subsection~\ref{s95subsecpflem}.

The passage to the limit under the infinite sums in the right-hand
side of \eqref{CH1_s95e500} is therefore justified. We have
\begin{equation}
\label{CH1_s95e600}
\begin{aligned}
    {}\lim_{n\ra\infty}&\sum_{ \genfrac{}{}{0pt}{}{\vka\in{}Q^m({}n/(2m))}
                     {\k_1+\cdots+\k_m=0} }\,\,
          \sum_{j=M_m(\vka)+1}^0
               \big|\,T_{{}n+j}B_{\vka}\,\big|\\
   %  &=\lim_{n\ra\infty}\sum_{ \genfrac{}{}{0pt}{}{\vka\in{}\Z^m}
    %                 {\k_1+\cdots+\k_m=0} }\,\,
     %     \sum_{j=M_m(\vka)+1}^0
      %         \,\Big|\,T_{{}n+j}\Big(B^{({}n/(2m))}\Big)_{\vka}\,\Big|,\\
     &=\sum_{ \genfrac{}{}{0pt}{}{\vka\in{}\Z^m}
                     {\k_1+\cdots+\k_m=0} }\,\,
          \sum_{j=M_m(\vka)+1}^0
               \,\lim_{n\ra\infty}
                 \Big|\,T_{{}n+j}\Big(B^{({}n/(2m))}\Big)_{\vka}\,\Big|,
\end{aligned}
\end{equation}
where the notation \eqref{s95eqABR} has been employed.
Due to the infinite smoothness of the full symbol $\sigma(B)$
and Egorov's theorem, and also by \eqref{s95eqL31}
$$
\begin{aligned}
       \lim_{n\ra\infty}
                 \Big|\,T_{{}n+j}\Big(B^{({}n/(2m))}\Big)_{\vka}\,\Big|
%      &
= \lim_{n\ra\infty}
                 \big|\,T_{n} B_{\vka}\,\big|%\\
%      &
= \Big|\,\res\big(A^{2-d}B_{\vka}\,\big)\,\Big|.
\end{aligned}
$$ 
Therefore the right-hand side of \eqref{CH1_s95e600} 
does not exceed 
$$
\begin{aligned}
      \sum_{ \vka\in\Z^m }\,\,
          \big(\,|\k_1|+\cdots+|\k_m|\,\big)
                  &\int_{S^*M} \Big|\,\sigma_{-d}\big(\,A^{2-d}
                       B_{\k_1}\cdots B_{\k_m}\,\big)\,\Big|\,dx\,d\xi\\
           &\leq c(B)\cdot m\cdot(C(B))^m,
\end{aligned}
$$
%for two certain constants depending only on $B$
again due to the smoothness of the symbol of any order 
of $B$,
%rapid decay of the Fourier coefficients $B_\k$,
%as $|\k|\ra\infty$, 
see subsection~\ref{s95subsecpflem}. 
We notice that $C(B)$ depends on
$ \int_{S^*M}\big|\nabla^3\sigma_{-d}(A^{2-d}B)\big|\,dxd\xi$,
$\|\nabla^3\sigma_0(B)\|_\infty$ 
and $\|\nabla^3\sub(B)\|_\infty$.
%$\|\nabla^3\sigma_{-2}(B)\|_\infty$.
\subsection{End of the proof of Theorem \ref{thm_main}}
%\label{CH1_subsecendofproof}
We have shown that there exist two constants
$c(B)$ and $C(B)$ (the latter constant is small together with 
the norm \eqref{eq_seminorm})
such that for all $m\in\N$, $m\geq2$, one has
\begin{equation}
\label{CH1_s95eR}
\begin{aligned}
    \Big|\tr&((P_n B P_n)^m) - \tr(P_n B^m P_n)%\\
%       &\qquad
- {}n^{d-1}\cdot\Upsilondva[z^m](B) 
                 -{}n^{d-2}\cdot \Upsilontri[z^m](B)\Big|\\
    &\quad\leq{}n^{d-3}\cdot c(B)\cdot m^4\cdot(C(B))^m,\quad n\ra\infty.
\end{aligned}
\end{equation}
Take now an arbitrary function
$f(z)=\sum_{m=1}^\infty{}c_mz^m\in\an$ which is analytic
on a neighborhood of $\{z\,:\,|z|<C(B)\}$.
Because the trace is a linear operation and the functionals 
$\Upsilon_l$, $l=2,3$, are linear in the first
argument, we can write
$$
\begin{aligned}
   \Big|\tr(\,f(P_n B P_n)\,) &- \tr(P_n f(B) P_n)%\\
%             &\qqquad
- {}n^{d-1}\cdot\Upsilondva[f(z)](B) 
                           -{}n^{d-2}\cdot \Upsilontri[f(z)](B)\Big|\\
         &= \bigg|\sum_{m=1}^\infty c_m
             \big[\tr(\,(P_n B P_n)^m\,) - \tr(P_n B^m P_n)\\
       &\qqquad
- {}n^{d-1}\cdot\Upsilondva[z^m](B) 
                  -{}n^{d-2}\cdot  \Upsilontri[z^m]B)\big]\bigg|.%\\
 %   &\leq\sum_{m=1}^\infty |c_m|\, 
 %
\end{aligned}
$$
Now we estimate the absolute value of the latter
sum by the sum of the absolute values of its terms
and employ \eqref{CH1_s95eR}. We obtain then
$$
\begin{aligned}
   \Big|\tr&(\,f(P_n B P_n)\,) - \tr(P_n f(B) P_n)%\\
%       &\qquad
- {}n^{d-1}\cdot\Upsilondva[f](B) -{}n^{d-2}\cdot
\Upsilontri[f](B)\Big|\\
    &\quad\leq{}n^{d-3}\cdot c(B)\sum_{m=1}^\infty
                    |c_m|\,m^4\,\cdot(C(B))^m,\quad n\ra\infty,
\end{aligned}
$$
 the series on the right is convergent due to the analyticity of
$f$ on a neighborhood of $\{z:|z|<C(B)\}$.
\subsection{Proof of Lemma \ref{CH1_s95l0}}
\label{s95subsecpflem}
We have to prove that for an arbitrary $B\in\Psi^0(M)$, $l=2,3$,
and any $L\in\N$ there exist two constants
$c_{l,L}(B)$ and $C_l(B)$ so that for $m\in\N$, $m\geq2$
$$
\begin{aligned}
    \Big|\,\Upsilon_l[z^m]\big(B^{(M)}\big) 
               &- \Upsilon_l[z^m](B)\,\Big|\\
        &\leq c_{l,L}(B)\cdot{}m^2\cdot(C_l(B))^m\cdot M^{-L},\quad M\ra\infty.
\end{aligned}
$$
We consider the case $l=2$ first. 
By the definition of $\Upsilon_2$ 
\begin{equation}
\label{eq_lem1}
\begin{aligned}
    \Big|\,\Upsilon_2[z^m]\big(B^{(M)}\big) 
               &- \Upsilon_2[z^m](B)\,\Big|\\
        &=  \Big|\, \sum_{\k_1+\cdots+\k_m=0}\,
            M_m(\vka)\,
   \int_{S^*M}\sigma_0\big(B_{\vka}^{(M)}-B_{\vka}\big)\,dxd\xi\,\Big|
\end{aligned}
\end{equation}
We have 
\begin{equation}
\label{eq_lem2}
\begin{aligned}  
      \null&B_{\k_1}^{(M)}\cdots{}B_{\k_m}^{(M)}
  -B_{\k_1}\cdots{}B_{\k_m}
    =\big(B_{\k_1}^{(M)}-B_{\k_1}\big)B_{\k_2}^{(M)}\cdots{}B_{\k_m}^{(M)}\\
       &+B_{\k_1}\big(B_{\k_2}^{(M)}-B_{\k_2}\big)B_{\k_3}^{(M)}\cdots{}B_{\k_m}^{(M)}+\cdots%\\
 %      &
+B_{\k_1}\cdots{}B_{\k_{m-1}}\big(B_{\k_m}^{(M)}-B_{\k_m}\big).
\end{aligned}
\end{equation}
Let us estimate the integral of the first summand, the other
$m-1$ are estimated analogously.
Note that 
$%\begin{equation}
%%\label{eq_lem3}
%\begin{aligned}
    \Big|\,M_m(\vka)\,\Big|\leq(1+|\k_1|)\cdots(1+|\k_m|)
%\end{aligned}
$.  %\end{equation}
Also because the sum of the indices in \eqref{eq_lem1}
is zero we have 
$$\begin{aligned}
     |\k_1|=|-\k_2-\cdots-\k_m|\leq|\k_2|+\cdots+|\k_m|
\end{aligned}
$$
and hence 
$$%\begin{equation}
%\label{eq_lem5}
%\begin{aligned}
    \Big|\,M_m(\vka)\,\Big|\leq(1+|\k_2|)^2\cdots(1+|\k_m|)^2.
%\end{aligned}
$$%\end{equation}
(We do not want $\k_1$ to appear on the right-hand side
in \eqref{eq_est1} below.)
Introduce the notation for $G\in\Psi^0(M)$
$$
      \|\sigma_0(G)\|_{\infty}:= \max_{(x,\xi)\in{}S^*M}
                \big|\sigma_0(G)(x,\xi)\big|.
$$
Hence 
\begin{equation}
\label{eq_est1}
\begin{aligned}
\Big|\, &\sum_{\k_1+\cdots+\k_m=0}\,
            M_m(\vka)\,
   \int_{S^*M}\sigma_0\big((B_{\k_1}^{(M)}-B_{\k_1})B_{\k_2}^{(M)}
                                \cdots{}B_{\k_m}^{(M)}\big)\,dxd\xi\,\Big|\\
   &\leq|S^*M|\sum_{\k_2,\cdots,\k_m\in\Z}\,
            (1+|\k_2|)^2\cdots(1+|\k_m|)^2\\
  &\times\big\|\sigma_0(B_{-\k_2­\cdots-\k_m}^{(M)}-B_{-\k_2­\cdots-\k_m})\big\|_{\infty}
      \|\sigma_0(B_{\k_2}^{(M)})\|_{\infty}\cdots{}\|\sigma_0(B_{\k_m}^{(M)})\|_{\infty}.
\end{aligned}
\end{equation}
For any $G\in\Psi^0(M)$ we write for brevity
$g_0:=\sigma_0(G)$, $g_{{\rm sub}}:=\sub(G)$,
introduce the notation 
\begin{equation}
\label{CH1_s5e2.5}
       (\ad A)(G):= [A,G]\in\Psi^0(M)
\end{equation}
and state the following important fact (see \cite[Lemma~1.3]{GO}).
\begin{lem}
\label{lem_rd}
 For any $G\in\Psi^0(M)$ and any $\nu\in\Z$
\begin{equation}
\label{s95eq030}
      \nu G_\nu = [A,G]_\nu,
\end{equation}
where the right-hand side is the $\nu$th Fourier coefficient 
of $[A,G]$. 
Also for any $L\in\N$ there exist finite $c_L(G)$ and $C_L(G)$ such that
$$%\begin{equation}
%\label{eq_rd1}
      \| \sigma_0(G_\nu)\|_{\infty} \leq C_L \nu^{-L},\qquad |\nu|\ra\infty,
$$%\end{equation}
and
\begin{equation}
\label{eq_rd1sum}
      \sum_{|\k|>M}\| \sigma_0\big(G_\k)\big)\|_{\infty} \leq C_L M^{-L},\qquad 
               M\ra\infty,
\end{equation}
and
$$%\begin{equation}
%\label{eq_rd1sub}
      \| \sub(G_\nu)\| \leq C_L \nu^{-L},\qquad |\nu|\ra\infty,
$$%\end{equation}
and
$$%\begin{equation}
%\label{eq_rd1ss}
      \sum_{|\k|> M}\| \sub\big(G_\k\big)\|_{\infty} 
           \leq C_L M^{-L},\qquad M\ra\infty,
$$%\end{equation}
and for the operator norm in $L^2(M)$
\begin{equation}
\label{eq_norm0}
      \sum_{|\k|> M} \|G_\k\|
           \leq C_L M^{-L},\qquad M\ra\infty.
\end{equation}
%($\| \cdot\|$ stands for the operator norm in $L^2(M)$).
\end{lem}
\begin{proof}
The equation \eqref{s95eq030} follows readily from the definition 
\eqref{CH1_s5e2.5}
and the fact that $\pi_\mu$ is the projector on the 
corresponding to the eigenvalue $\mu\in\N$ eigenspace 
of the operator $A$. %The other two equations, see \cite{GO}.
Due to the infinite smoothness of
$\sigma_0(G)$ we can repeat the operation \eqref{s95eq030}
any finite number of times, and consequently 
the norm $\|\sigma(G_\nu)\|_{\infty}$
decays rapidly, as $|\nu|\ra\infty$.
By the fast decay of $\|\sigma_0(G_\k)\|_\infty$, as $|\k|\ra\infty$, 
for any $L\in\N$ there is a constant $C_L(G)$ such that 
for all $(x,\xi)\in{}S^*M$ and $0\leq{}t\leq2\pi$
$$%\begin{equation}
%\label{s95eq020}
\begin{aligned}
         \big|(g^{(M)})_0^{t_1}-g_0^{t_1}\big| 
           &= \Big| \sum_{|\k|>M} \sigma_0(G_\k)(\Theta^t(x,\xi)\Big|\\
           &\leq  \sum_{|\k|>M} \|\sigma_0(G_\k)\|_\infty 
          \leq     C_L(G) \cdot M^{-L},\quad M\ra\infty
\end{aligned}
$$%\end{equation}
which together with 
Lemma \ref{lem0sub} proves \eqref{eq_rd1sum}. 
The result for the subprincipal symbol
also follows from its infinite smoothness and Lemma \ref{lem0sub}.
For \eqref{eq_norm0}, see \cite[Lemma~1.3]{GO}.
\end{proof}
For any $G\in\Psi^0(M)$ 
we now define inductively $(\ad{}A)^N(G)$, $N\in\N$,
and note that 
$$%\begin{equation}
%\label{s95eqadn}
     (\ad{}A)^N(G_\nu) = ((\ad{}A)^N(G))_\nu,\quad \nu\in\Z,\,N\in\N.
$$%\end{equation}
For $N=1$ this follows readily from \eqref{CH1_s5e2.5},
% follows easily as above that
%$[A,G_\nu]=([A,G])_\nu$, 
for $N\geq2$ we proceed by induction. 
\begin{lem}
For any $H\in\Psi^0(M)$ and $\nu\in\Z$
$$%\begin{equation}
%\label{s95eqtri}
            \|\sigma_0(H_\nu)\|_\infty \leq \|\sigma_0(H)\|_\infty
$$%\end{equation}
and
$$%\begin{equation}
%\label{s95eqtrip}
            \|\sub(H_\nu)\|_\infty \leq \|\sub(H)\|_\infty
$$%\end{equation}
\end{lem}
\begin{proof}
 Follows from Egorov's theorem \eqref{eqEgo} 
and \eqref{eqEsub}.
\end{proof}
\begin{lem}
\label{lem2.9}
For any $G\in\Psi^0(M)$
\begin{equation}
\label{s95eqPP}
       \sum_{\nu\in\Z}\big\| \sigma_0(G_\nu)\big\|_\infty 
           \leq \|\sigma_0(G_0)\|_\infty 
     + \frac{\pi^2}3\,\big\| \sigma_0((\ad{}A)^2(G))\big\|_\infty
\end{equation}
and 
$$%\begin{equation}
%\label{s95eqPPsub}
       \sum_{\nu\in\Z}\big\| \sub(G_\nu)\big\|_\infty 
           \leq \|\sub(G_0)\|_\infty 
     + \frac{\pi^2}3\,\big\| \sub((\ad{}A)^2(G))\big\|_\infty
$$%\end{equation}
and for the operator norm in $L^2(M)$
\begin{equation}
\label{eq_norm}
       \sum_{\nu\in\Z}\big\| G_\nu\big\|
           \leq \|G_0\|
     + \frac{\pi^2}3\,\big\| (\ad{}A)^2(G)\big\|.
\end{equation}
\end{lem}
\begin{proof}
Applying \eqref{s95eq030} twice we get
$$%\begin{equation}
%\label{s95eqraz}
   G_\nu = \nu^{-2}\,[A,[A,G]]_\nu,\quad\nu\neq0,
$$%\end{equation}
which implies \eqref{s95eqPP}
because $\sum_{\nu\neq0}\nu^{-2}=\pi^2/3$.
For \eqref{eq_norm}, see \cite[Lemma~1.3]{GO}.
\end{proof}

Now we estimate \eqref{eq_est1} by
$$%\begin{equation}
%\label{eq_est2}
\begin{aligned}
   \sum_{\k_1\in\Z}&\big\|\sigma_0((B^{(M)}-B)_{\k_1})\big\|_{\infty}\\
   &\times\bigg(\sum_{\k_2\in\Z}(1+|\k_2|)^2\|\sigma_0(B_{\k_2}^{(M)})\|_{\infty}\bigg)
     \cdots\bigg(\sum_{\k_m\in\Z}\,(1+|\k_m|)^2\|\sigma_0(B_{\k_m}^{(M)})\|_{\infty}\bigg).
\end{aligned}
$$%\end{equation}
By \eqref{eq_rd1sum} for any $L\in\N$ there is $c_{2,L}(B)$
such that, as $M\ra\infty$,
$$%\begin{equation}
%\label{eq_est3}
%\begin{aligned}
   \sum_{\k_1\in\Z}\big\|\sigma_0((B^{(M)}-B)_{\k_1})\big\|_{\infty}
   =\sum_{|\k_1|>M}\big\|\sigma_0(B_{\k_1})\big\|_{\infty}\leq c_{2,L}(B)\,M^{-L}.
%\end{aligned}
$$%\end{equation}
Also
$$%\begin{equation}
%\label{eq_est5}
\begin{aligned}
   \sum_{\k_2\in\Z}(1+|\k_2|)^2\|\sigma_0(B_{\k_2}^{(M)})\|_{\infty}
     =\sum_{\k_2\in\Z}(1+2|\k_2|+|\k_2|^2)\|\sigma_0(B_{\k_2}^{(M)})\|_{\infty}
\end{aligned}
$$%\end{equation}
and for example
\begin{equation}
\label{eq_est5p}
\begin{aligned}
   \sum_{\k_2\in\Z}|\k_2|^2\|\sigma_0(B_{\k_2}^{(M)})\|_{\infty}
         &=\sum_{|\k_2|\leq{}M}\|\sigma_0(\k_2^2B_{\k_2})\|_{\infty}\\
         &=\sum_{|\k_2|\leq{}M}\|\sigma_0((\ad A)^2(B))_{\k_2})\|_{\infty}\\
         &\leq\sum_{\k_2\in\Z}\|\sigma_0((\ad A)^2(B))_{\k_2})\|_{\infty}=:C(B),
\end{aligned}
\end{equation}
where $C(B)$ is estimated by $\|\nabla^{4}\sigma_0(B)\|_\infty$, 
in view of Lemma~\ref{lem2.9}.
We have proved that 
$$%\begin{equation}
%\label{eq_lem1rim}
\begin{aligned}
    \Big|\,\Upsilon_2[z^m]\big(B^{(M)}\big) 
               &- \Upsilon_2[z^m](B)\,\Big|
        \leq c_{2,L}(B)\cdot m\cdot(C_2(B))^m\cdot M^{-L},\qquad M\geq1. 
\end{aligned}
$$%\end{equation}
where $C_2(B)$ depends on $\|\nabla^{4}\sigma_0(B)\|_\infty$.

Now the proof for $\Upsilon_3$. In view of \eqref{CH1_s9e24}, the 
functional $\Upsilon_3$,
which is the coefficient of $n^{d-2}$ in \eqref{eq_together}, 
contains three parts.
The first one is up to a constant as the right-hand side of \eqref{eq_lem1},
and is estimated in the same way as above.
The second part is 
$$%\begin{equation}
%\label{eq_lem2pp}
%\begin{aligned}
 %   \Big|\,\Upsilon_2[z^m]\big(B^{(M)}\big) 
  %             &- \Upsilon_2[z^m](B)\,\Big|\\
%        &=  
\Big|\, \sum_{\k_1+\cdots+\k_m=0}\,
            (M_m(\vka))^2\,
   \int_{S^*M}\sigma_0\big(B_{\vka}^{(M)}-B_{\vka}\big)\,dxd\xi\,\Big|
%\end{aligned}
$$%\end{equation}
and the only difference with the above argument
is that here we have $(1+|\k_l|)^4$, $l=2,\cdots,m$ in the estimates,
and consequently 
$$
       C(B):= \sum_{\k\in\Z}\|\sigma_0\big( ((\ad A)^4(B))_{\k}\big)\|_{\infty} 
$$
can be estimated in terms of $\|\nabla^6\sigma_0(B)\|_\infty$, 
by Lemma \ref{lem2.9}.
The absolute value of the third part of
$
  \Upsilon_3[z^m]\big(B^{(M)}\big) -\Upsilon_3[z^m](B)
$
is estimated by
$$%\begin{equation}
%\label{eq_lem7}
\begin{aligned}
  \Big|\, \sum_{\k_1+\cdots+\k_m=0}\,
            M_m(\vka)\,
   \int_{S^*M}\sub\big(B_{\vka}^{(M)}-B_{\vka}\big)\,dxd\xi\,\Big|
\end{aligned}
$$%\end{equation}
We use \eqref{eq_lem2} and estimate for instance the
integral of the first of the $m$ summands
\begin{equation}
\label{eq_lem2p}
\begin{aligned}  
    \sub\big(\,(B_{\k_1}^{(M)}&-B_{\k_1})B_{\k_2}^{(M)}\cdots{}B_{\k_m}^{(M)}\,\big)
        = \sub(B_{\k_1}^{(M)}-B_{\k_1})\,
               \sigma_0\big(B_{\k_2}^{(M)}\cdots{}B_{\k_m}^{(M)}\big)\\
       &+ \sigma_0(B_{\k_1}^{(M)}-B_{\k_1})\,
               \sub\big(B_{\k_2}^{(M)}\cdots{}B_{\k_m}^{(M)}\big)\\
       &+ \frac1{2i}\,\big\{\sigma_0(B_{\k_1}^{(M)}-B_{\k_1}),\,
               \sigma_0(B_{\k_2}^{(M)}\cdots{}B_{\k_m}^{(M)})\big\}
\end{aligned}
\end{equation}
Denote $G:=(\ad{}A)^2(B)$. By Lemma \ref{lem2.9} 
for any $L\in\N$ %analogously to \eqref{eq_est5p} 
we obtain
$$%\begin{equation}
%\label{eq_lem10}
\begin{aligned}
    \Big|\,\sum_{\k_1+\cdots+\k_m=0}\,
            M_m(\vka)\,
   \int_{S^*M}&\sub(B_{\k_1}^{(M)}-B_{\k_1})\,
               \sigma_0\big(B_{\k_2}^{(M)}\cdots{}B_{\k_m}^{(M)}\big)\,dxd\xi\,\Big|\\
        &\leq c_L\cdot m\cdot(C(G))^m\cdot M^{-L},\qquad M\geq1. 
\end{aligned}
$$%\end{equation}
We note that $C(G(B))$ can be estimated in terms 
of $\|\nabla^2\sub(B)\|_\infty$ and $\|\nabla^4\sigma_0(B)\|_\infty$.
Now for the Poisson bracket term in \eqref{eq_lem2p}.
\begin{lem}
\label{lem2.10}
For any $G\in\Psi^0(M)$ and $L\in\N$
there exist $c_L(G)$ and $C_L(G)$ such that
$$%\begin{equation}
%\label{eq_pr1}
      \| \partial_x\sigma_0(G_\nu)\|_{\infty} \leq C_L \nu^{-L},\qquad |\nu|\ra\infty
$$%\end{equation}
and
$$%\begin{equation}
%\label{eq_pr2}
      \| \partial_\xi\sigma_0(G_\nu)\|_{\infty} \leq C_L \nu^{-L},\qquad |\nu|\ra\infty
$$%\end{equation}
and
$$%\begin{equation}
%\label{eq_lem15}
        C_1(G):= \sum_{\k\in\Z}   \|\partial_x\sigma_0(G_\k)\|_\infty <\infty
$$%\end{equation}
and
$$%\begin{equation}
%\label{eq_lem16}
           C_2(G):= \sum_{\k\in\Z}   \|\partial_\xi\sigma_0(G_\k)\|_\infty <\infty
$$%\end{equation}
\end{lem}
Now we recall \eqref{CH1_s9e10}, estimate the absolute
value of the Poisson brackets 
by the products of the supremum norms
of the derivatives, employ Lemma~\ref{lem2.10} and conclude
that the contribution of the Poisson brackets
involves $(C(B))^m$, where $C(B)$ depends on $\|\nabla^5\sigma_0(B)\|_\infty$. 
Thus for a certain $C_3(B)<\infty$ and any $L\in\N$
$$%\begin{equation}
%\label{eq_lem1lst}
\begin{aligned}
    \Big|\,\Upsilon_3[z^m]\big(B^{(M)}\big) 
               &- \Upsilon_3[z^m](B)\,\Big|
        \leq c_{3,L}(B)\cdot m^2\cdot(C_3(B))^m\cdot M^{-L},\qquad M\geq1, 
\end{aligned}
$$%\end{equation}
where $C_3(B)$ depends on $\|\nabla^{6}\sigma_0(B)\|_\infty$ and
$\|\nabla^2\sub(B)\|_\infty$.
The factor $m^2$ is due to the fact that there are $m$
terms in the expansion of the Poisson bracket in \eqref{eq_lem2p}.
\section%[How the basic functors $\tilde\Phi_j$, $j\in\N$, do arise]
{Contribution of %the principal symbol of $B_{\k_1}\cdots{}B_{\k_m}$}
$\sigma_0\big(B_{\k_1}\cdots{}B_{\k_m}\big)$}
%{An auxiliary calculation: how the basic functors $\tilde\Phi_j$, $j\in\N$,
%arise from an application of gHD in an abstract setting,
%and two applications}
\label{CH1_s6aux}
Let $B\in\Psi^0(M)$ (we redenote $G$ by $B$).
Our goal here is to calculate 
\begin{equation}
\label{CH1_s6e3}
%\begin{aligned}
   \sum_{\k_1+\cdots+\k_m=0}\,
            (\,M_m({\overline{\k}})\,)^n \,
   \int_{S^*M}\sigma_0(B_{{\overline{\k}}})\,dxd\xi
%\end{aligned}
\end{equation}
for $n=1,2$. Assume for the moment that $n$ is any natural number.
Recall that $S_m$ denotes the set of all permutations $\tau$
of the numbers $1,\cdots,m$. Denote
$    {\overline{\k}}_\tau  := (\k_{\tau_1},\ldots,\k_{\tau_m}).
$
Note that both the second factor in \eqref{CH1_s6e3} 
and the domain of summation are symmetric in $\k$.
Therefore
$$%\begin{equation}
%\label{CH1_s6e5}
\begin{aligned}
     \null&\sum_{\k_1+\cdots+\k_m=0}\,\frac1{m!}\sum_{\tau\in{}S_m}
            (\,M_m({\overline{\k}_\tau})\,)^n \,\int_{S^*M}\sigma_0(B_{{\overline{\k}}})\,dxd\xi\\
    &= \sum_{\k_1+\cdots+\k_m=0}\,\frac1{m!}\sum_{\tau\in{}S_m}
         \sum_{p=1}^{m-1}
           \Big[ (\,M_{p}({\overline{\k}_\tau})\,)^n 
                  - (\,M_{p-1}({\overline{\k}_\tau})\,)^n \Big]
                       \,\int_{S^*M}\sigma_0(B_{{\overline{\k}}})\,dxd\xi.
\end{aligned}
$$%\end{equation}
Now an application of the gHD (Theorem \ref{CH1_S2TH1}) gives
$$%\begin{equation}
%\label{CH1_s6e7}
\begin{aligned}
  \sum_{\k_1+\cdots+\k_m=0}
         \frac1{m!} &\sum_{\tau\in{}S_m}
          \sum_{p=1}^{m-1}
    \sum_{j=1}^{\min(p,n)}
          \,\frac1{j!}\sum_{   \genfrac{}{}{0pt}{}  {  k_1,\cdots,k_j \geq 1 }
                  { k_1+\cdots+k_j = p }
           } \,
  \sum_{    \genfrac{}{}{0pt}{} {  l_1,\cdots,l_j \geq 1 }
                  { l_1+\cdots+l_j = n }
           } \,
  \lp \genfrac{}{}{0pt}{}{n}{ l_1,\cdots, l_j} \rp \\
      &\times\frac{(-(k_1({\overline{\k}}_\tau))_-)^{l_1}}{k_1}\cdot\ldots\cdot
                \frac{(-(k_j({\overline{\k}}_\tau))_-)^{l_j}}{k_j}
               \cdot\int_{S^*M}\sigma_0(B_{{\overline{\k}}})\,dxd\xi,
\end{aligned}
$$%\end{equation}
where the notation \eqref{CH1_s1e6prim} has been used. 
We make the sum $\frac1{m!}\sum_{\tau\in{}S_m}$ the outmost one
and then drop it, because it is the number of $\k$'s that matters, not their
indices (symmetricity of the last factor and of the domain of summation
again).  We get 
\begin{equation}
\label{CH1_s6e8}
\begin{aligned}
  \sum_{\k_1+\cdots+\k_m=0}
   &\sum_{p=1}^{m-1}
    \sum_{j=1}^{\min(p,n)}
          \,\frac1{j!} \sum_{     \genfrac{}{}{0pt}{}{  k_1,\cdots,k_j \geq 1 }
                  { k_1+\cdots+k_j = p }
           } \,
  \sum_{    \genfrac{}{}{0pt}{} {  l_1,\cdots,l_j \geq 1 }
                  { l_1+\cdots+l_j = n }
           } \,
  \lp \genfrac{}{}{0pt}{}{n}{ l_1,\cdots, l_j} \rp \\
      &\times\frac{(-(k_1({\overline{\k}}))_-)^{l_1}}{k_1}\cdot\ldots\cdot
                \frac{(-(k_j({\overline{\k}}))_-)^{l_j}}{k_j}
               \cdot\int_{S^*M}\sigma_0(B_{{\overline{\k}}})\,dxd\xi.
\end{aligned}
\end{equation}
We have $m-1$ independent summation variables $\overline{\k}$.
Let us make the change of variables
$$%\begin{equation}
%\label{CH1_s6e9}
\begin{aligned}
  \nu_1 &:= k_1({\overline{\k}}) = \k_1+\cdots+\k_{k_1}\\
  \nu_2 &:= k_2({\overline{\k}}) = \k_{k_1+1}+\cdots+\k_{k_1+k_2}\\
            &\cdots\\
   \nu_j &:= k_j({\overline{\k}}) = \k_{k_1+\cdots+k_{j-1}+1}
                                       +\cdots+\k_{k_1+\cdots+k_{j-1}+k_j}
\end{aligned}
$$%\end{equation}
for $j=1,\cdots,\min(p,n)$. We also interchange the summations 
over $p$ and $j$,
introduce the new variable $k_{j+1}:=m-p$,
and perform the $m-j-1$ unconditional summations in \eqref{CH1_s6e8}.
We need the following important fact \cite[after Lemma~1.2]{GO}
which states that 
the convolution on the Fourier transform side
corresponds to raising $B$ to a power on the original side.
\begin{lem}
\label{s6auxlemGO} 
For any $B\in\Psi^0(M)$ and arbitrary $j\in\N$ and $\nu\in\Z$ one has
$$%\begin{equation}
%\label{CH1_s6e100}
    \sum_{\k_1+\cdots+\k_j=\nu}  B_{\k_1}\cdots B_{\k_j}
           = (B^j)_\nu.
$$%\end{equation}
\end{lem}
Notice also that by Egorov's theorem \eqref{eqEgo}
and because $\sigma_0(B^k)=(\sigma_0(B))^k$,
$$
           \sigma_0((B^k)_\nu)(x,\xi) = \int_{0}^{2\pi} (b_0^{t}(x,\xi))^k\,e^{it\nu}\,\frac{dt}{2\pi}
          = \int_{0}^{2\pi} (b_0^{k}(x,\xi))^t\,e^{it\nu}\,\frac{dt}{2\pi}.
$$
Then \eqref{CH1_s6e8} becomes
\begin{equation}
\label{CH1_s6e17}
\begin{aligned}
   \null&\sum_{j=1}^{\min(n,m-1)}
          \frac1{j!} \sum_{\nu_1,\cdots,\nu_j}
         \sum_{     \genfrac{}{}{0pt}{}{  l_1,\cdots,l_j \geq 1 }
                  { l_1+\cdots+l_j = n }
           } \,
                \lp \genfrac{}{}{0pt}{}{n}{ l_1,\cdots, l_j} \rp
             \cdot(-(\nu_1)_-)^{l_1}\cdot\ldots\cdot(-(\nu_j)_-)^{l_j}\\
     &\times\int_{S^*M} dxd\xi\,
   \int_0^{2\pi}\!\cdots\!\int_0^{2\pi}
                 e^{i(\nu_1t_1+\cdots+\nu_jt_j
                              -(\nu_1+\cdots+\nu_j)t_{j+1})}\,
   \frac{dt_1}{2\pi}\cdots
               \frac{dt_j}{2\pi}\frac{dt_{j+1}}{2\pi}\\
        &\times F_j[z^m](b_0^{t_1},\cdots,b_0^{t_j},b_0^{t_{j+1}}\,)
\end{aligned}
\end{equation}
where we have introduced a linear 
$(j+1)$-map $F_{j+1}$, whose
action on $f(z)=z^m$, $m\geq{}j+1$,
is prescribed by 
\begin{equation}
\label{CH1_s6e16}
%\begin{aligned}
  F_{j+1}[z^m](x_1,\cdots,x_j,x_{j+1}) %&
           :=\sum_{     \genfrac{}{}{0pt}{}{k_1,\cdots,k_j,k_{j+1}\geq1}
                  { k_1+\cdots+k_j+k_{j+1} = m }
           } \frac{x_1^{k_1}}{k_1}\cdot\ldots\cdot
                      \frac{x_j^{k_j}}{k_j}\cdot
                        x_{j+1}^{k_{j+1}}
%\end{aligned}
\end{equation}
and $F_{j+1}[z^m]:=0$ for $m=1,\cdots,j$.
\begin{rem} 
It is important that for any fixed $n$ the sum
over $j$ in \eqref{CH1_s6e17} termina\-tes at $n$, no matter how large $m$ is.
\end{rem}
We need the following fact \cite[after~(1.3)]{GO}.
\begin{lem}
\label{lemsymm}
Let $M$ be a Zoll manifold and $A$ defined by \eqref{eq_A}.
Assume $\sigma_1(A)(x,\xi)=\sigma_1(A)(x,-\xi)$ 
for all $(x,\xi)\in{}T^*M$.
Let $f,g\in{}C^\infty(S^*M)$, $\nu\in\Z$,
and recall the notation \eqref{eq_FcfZ}.
Then
$$
    \int_{S^*M} \hat{f}_\nu\, \hat{g}_{-\nu}\, dxd\xi 
   = \int_{S^*M} \hat{f}_{-\nu}\, \hat{g}_\nu\, dxd\xi.
$$
\end{lem}
We need to compute \eqref{CH1_s6e17} only for $n=1,2$.
The computation for $n=1$ has been carried out in \cite{GO}.
We repeat it here for the sake of completeness.
For $n=1$ \eqref{CH1_s6e17} becomes
\begin{equation}
\label{CH1_s6e17prim}
\begin{aligned}
   \sum_{\nu_1}
         &(-(\nu_1)_-)\int_{S^*M} dxd\xi\,
   \int_0^{2\pi}\int_0^{2\pi}
                 e^{i\nu_1(t_1-t_2)}\,
            \sum_{    \genfrac{}{}{0pt}{} {k_1,k_{2}\geq1}
                  { k_1+k_{2} = m }
           } \frac{(b_0^{k_1})^{t_1}}{k_1}(b_0^{k_2})^{t_2}
   \frac{dt_1}{2\pi}\frac{dt_2}{2\pi}\\
   &=\sum_{\nu_1}
         (-(\nu_1)_-)\int_{S^*M} dxd\xi\,
   \sum_{    \genfrac{}{}{0pt}{} {k_1,k_{2}\geq1}
                  { k_1+k_{2} = m }
           } \bigg(\frac{\widehat{{b_0^{k_1}}}}{k_1}\bigg)_{\nu_1}\,
     (\widehat{{b_0^{k_2}}})_{-\nu_1}
\end{aligned}
\end{equation}
This can be rewritten as 
\begin{equation}
\label{CH1_s6e17bis}
%\begin{aligned}
\sum_{\nu}
        (-(\nu)_-)\int_{S^*M} dxd\xi\,
   \sum_{ j=1}^{m-1}
   \bigg(\frac{\widehat{{b_0^{j}}}}{j}\bigg)_{\nu}\,
     (\widehat{{b_0^{m-j}}})_{-\nu}
%\end{aligned}
\end{equation}
which by Lemma  \ref{lemsymm} equals
\begin{equation}
\label{CH1_s6e17t}
\begin{aligned}
   \sum_{\nu}
        &(-(\nu)_-)\int_{S^*M} dxd\xi\,
  \sum_{ j=1}^{m-1}
      (\widehat{{b_0^{m-j}}})_{\nu}\,\bigg(\frac{\widehat{b_0^{j}}}{j}\bigg)_{-\nu}\\
   &=\sum_{\nu}
        (-(\nu)_-)\int_{S^*M} dxd\xi\,
  \sum_{ j=1}^{m-1}
      (\widehat{{b_0^{j}}})_{\nu}\,\bigg(\frac{\widehat{b_0^{m-j}}}{m-j}\bigg)_{-\nu}
\end{aligned}
\end{equation}
where the last equality is obtained after the change of summation index $j\ra{m-j}$.
Now we sum \eqref{CH1_s6e17bis} and \eqref{CH1_s6e17t}
and use $\frac1j+\frac1{m-j}=\frac{m}{j(m-j)}$ 
and Lemma~\ref{lemsymm} to conclude
that \eqref{CH1_s6e17prim} equals
$$%\begin{equation}
%\label{CH1_s6e17f}
\begin{aligned}
 -\frac{m}2\sum_{\nu}
        &(-(\nu)_-)\int_{S^*M} dxd\xi\,
  \sum_{ j=1}^{m-1}
      \bigg(\frac{\widehat{b_0^{j}}}{j}\bigg)_{\nu}\,
            \bigg(\frac{\widehat{b_0^{m-j}}}{m-j}\bigg)_{-\nu}\\
   &=\sum_{k=1}^\infty k\,
        \int_{S^*M} dxd\xi\,
          \int_0^{2\pi}\int_0^{2\pi} e^{ik(t_1-t_2)}\, {W}_2[z^m](b_0^{t_1},b_0^{t_2})\,
               \frac{dt_1}{2\pi}\frac{dt_2}{2\pi},
\end{aligned}
$$%\end{equation}
where following \cite{LRS} we have noticed that $W_2$ 
defined by \eqref{CH1_s5ieW} satisfies
$$
    W_2[z^m](x_1,x_2) = \frac{m}2 \sum_{j=1}^{m-1} \frac{x_1^{j}}{j}
      \,\frac{x_2^{m-j}}{m-j}.
$$

Now for the case when $n=2$ in \eqref{CH1_s6e17}.
In this case $j$ can be either $1$ or $2$. In the first case
we repeat the above argument, the corresponding contribution becomes
$$%\begin{equation}
%\label{CH1_s6e17s}
\begin{aligned}
   \sum_{k=1}^\infty k^2\,
        \int_{S^*M} dxd\xi\,
          \int_0^{2\pi}\int_0^{2\pi} e^{ik(t_1-t_2)}\, {W}_2[z^m](b_0^{t_1},b_0^{t_2})\,
               \frac{dt_1}{2\pi}\frac{dt_2}{2\pi}.
\end{aligned}
$$%\end{equation}
For $j=2$, \eqref{CH1_s6e17} becomes
\begin{equation}
\label{CH1_s6e17se}
\begin{aligned}
   \sum_{\nu_1,\nu_2}
         &(-(\nu_1)_-)\,(-(\nu_2)_-)\int_{S^*M} dxd\xi\,\\
   &\times\int_0^{2\pi}\int_0^{2\pi}\int_0^{2\pi}
                 e^{i(\nu_1t_1+\nu_2t_2
                              -(\nu_1+\nu_2)t_{3})} \, 
    F_3[z^m](b_0^{t_1},b_0^{t_2},b_0^{t_{3}})
   \frac{dt_1}{2\pi}\frac{dt_2}{2\pi}\frac{dt_{3}}{2\pi},
\end{aligned}
\end{equation}
where by \eqref{CH1_s6e16}
\begin{equation}
\label{CH1_s6e16biss}
%\begin{aligned}
  F_{3}[z^m](x_1,x_2,x_3) %&
           :=\sum_{     \genfrac{}{}{0pt}{}{k_1,k_2,k_{3}\geq1}
                  { k_1+k_2+k_{3} = m }
           } \frac{x_1^{k_1}}{k_1}\,
                      \frac{x_2^{k_2}}{k_2}\cdot
                        x_{3}^{k_{3}},\qquad m\geq3,
%\end{aligned}
\end{equation}
and $F_{3}[z^m]:=0$ for $m=1,2$.
Let us redenote $W_3:=F_3$.
The formula \eqref{CH1_s5ieFprin} for $W_3[f]$,
$f\in\an$, is justified in Section~\ref{CH1_s6Phi}. 
Now the last detail, the following can be proved 
in the same way as Lemma~\ref{lemsymm}.
For any $f,g,h\in{}C^\infty(S^*M)$ and $\lambda,\mu,\nu\in\Z$
\begin{equation}
\label{eq_addit}
      \int_{S^*M} \hat{f}_\lambda\,\hat{g}_\mu\,\hat{h}_\nu\,dxd\xi
            =\int_{S^*M} \hat{f}_{-\lambda}\,\hat{g}_{-\mu}\,\hat{h}_{-\nu}\,dxd\xi.
\end{equation}
In view of \eqref{eq_addit} we rewrite \eqref{CH1_s6e17se} as
$$%\begin{equation}
%\label{CH1_s6e17sw}
%\begin{aligned}
   \sum_{k,l=1}^\infty k\,l
         \int_{S^*M} dxd\xi\,\int_0^{2\pi}\!\int_0^{2\pi}\!\int_0^{2\pi}
                 e^{i(kt_1+lt_2-(k+l)t_3)} \, 
    W_3[z^m](b_0^{t_1},b_0^{t_2},b_0^{t_3})
   \frac{dt_1}{2\pi}\frac{dt_2}{2\pi}\frac{d_3}{2\pi}.
%\end{aligned}
$$%\end{equation}
\section
%[Contribution of the symmetric part of the subprincipal symbol]
{Contribution of the symmetric part of 
$\sub\!\big(B_{\k_1}\cdots{}B_{\k_m}\big)$}
%the subprincipal symbol of $B_{\k_1}\cdots{}B_{\k_m}$}
%{An auxiliary calculation for the third order Szeg\"o functional
%corresponding to the symmetric part of the subprincipal symbol}
\label{CH1_s6sub}
Let us choose any $B\in\Psi^0(M)$.
The goal of this section is to calculate 
for an arbitrary $m\in\N$, $m\geq2$,
the sum
\begin{equation}
\label{CH1_s6sube3}
\begin{aligned}
   \sum_{\k_1+\cdots+\k_m=0}\,
           &M_m({\overline{\k}})\,
         \int_{S^*M}dxd\xi\,\\
         &\times\sum_{r=1}^m \sigma_0(B_{\k_1}\cdots{}B_{\k_{r-1}})
         \,\sub\big(B_{\k_r}\big)
      %\big[\sigma_1(A)\,\sub(B_{\k_r})+\frac{d-1}{2i}\{\sigma_1(A),\sigma_0(B_{\k_r})\}\big]
          \,\sigma_0(B_{\k_{r+1}}\cdots{}B_{\k_m}),
\end{aligned}
\end{equation}
As in Section \ref{CH1_s6aux}, the key observation is 
that both the second factor in 
\eqref{CH1_s6sube3} and the domain of summation are symmetric
in $\overline{\k}$. We 
permute all the $\k$'s in the first factor in \eqref{CH1_s6e3},  make
use of Theorem \ref{CH1_S2TH1} for the power $n=1$,
that is the classic HD,
$$%\begin{equation}
%\label{CH1_s6sube1}
\begin{aligned}
 \sum_{\tau\in{}S_m} M_m(\k_\tau) 
     = \sum_{\tau\in{}S_m} \sum_{j=1}^m
       \frac{-(\k_{\tau_1}+\cdots+\k_{\tau_j})_-}{j},
\end{aligned}
$$%\end{equation}
and drop the summation over $S_m$, because, as in Section \ref{CH1_s6aux}, 
it is the number of $\k$'s in a particular sum that matters, 
not their indices (again, the domain of summation and the second factor are still
symmetric). 
After that we get
$$%\begin{equation}
%\label{CH1_s6sube9}
\begin{aligned}
     \null\int_{S^*M} dxd\xi\,
          &\sum_{j=1}^{m-1}\sum_{r=1}^{m}\sum_{\k_1+\cdots+\k_m=0}
                   \frac{-(\k_{1}+\cdots+\k_{j})_-}{j} \\
       &\times\sigma_0(B_{\k_1}\cdots{}B_{\k_{r-1}})
         \,\sub\big(B_{\k_r}\big)
          \,\sigma_0(B_{\k_{r+1}}\cdots{}B_{\k_m}).
\end{aligned}
$$%\end{equation}
We split the above sum into two: for $r=1,\cdots,j$
and for $r=j+1,\cdots,m$. In the first case we set 
$$
            \nu:=\k_1+\cdots+\k_j,\quad -\mu+\nu:=\k_r.
$$
Then 
$$
         \k_1+\cdots+\k_{r-1}+\k_{r+1}+\cdots+\k_j = \mu,\quad 
               \k_{j+1}+\cdots+\k_m=-\nu.
$$
We carry out the independent summations
and recall Egorov's theorem
\eqref{eqEgo} and Lemma \ref{s6auxlemGO} to write the sum over $r=1,\cdots,j$ as
\begin{equation}
\label{CH1_s6sube9p}
\begin{aligned}
      \int_{S^*M} dxd\xi\,
          \sum_{j=1}^{m-1}\,\frac1j\sum_{r=1}^{j}\sum_{\nu,\mu}(-(\nu)_-) %\\
       \sigma_0((B^{j-1})_{\nu-\mu})
               \sub(B_{\mu})\sigma_0((B^{m-j})_{-\nu})
\end{aligned}
\end{equation}
Now let us write $h_0^t:=\sigma_0(H)(\Theta^t(x,\xi))$
and $h_{{\rm sub}}^t:=\sub(H)(\Theta^t(x,\xi))$, $H\in\Psi^0(M)$.
We have
$$%\begin{equation}
%\label{eq_sub1}
\begin{aligned}
     \sum_{\mu}
       \sigma_0((&B^{j-1})_{\nu-\mu})\,
               \sub(B_{\mu})=\sum_{\mu}\int_0^{2\pi}\int_0^{2\pi}
           e^{is(\nu-\mu)+it\mu} (b_0^s)^{j-1}\,b_{{\rm sub}}^t\,
                     \frac{ds}{2\pi}\frac{dt}{2\pi}\\
      &=\int_0^{2\pi}
           e^{is\nu} (b_0^s)^{j-1}\,b_{{\rm sub}}^s\,
                     \frac{ds}{2\pi}%\\
%    &
= \frac1j\sub\big( (B^j)_\nu\big).
\end{aligned}
$$%\end{equation}
Then \eqref{CH1_s6sube9p} becomes
\begin{equation}
\label{eq_100}
\begin{aligned}
      \int_{S^*M} dxd\xi\,
          \sum_{j=1}^{m-1}\,\frac1j\sum_{\nu}(-(\nu)_-) %\\
      \sub\big( (B^j)_\nu\big)\,\sigma_0((B^{m-j})_{-\nu})
\end{aligned}
\end{equation}
Analogously, in the case $r=j+1,\cdots,m$ 
we set 
$$
            \nu:=\k_1+\cdots+\k_j,\quad -\mu-\nu:=\k_r.
$$
Then 
$$
         \k_{j+1}+\cdots+\k_{r-1}+\k_{r+1}+\cdots+\k_j = \mu,\quad 
               \k_{j+1}+\cdots+\k_m=-\nu.
$$
Then the sum for $r=j+1,\cdots,m$ becomes
$$%\begin{equation}
%\label{CH1_s6sube9pp}
\begin{aligned}
      \int_{S^*M} dxd\xi\,
          &\sum_{j=1}^{m-1}\,\frac1j\sum_{r={j+1}}^{m}\sum_{\nu,\mu}(-(\nu)_-) 
       \sigma_0((B^{j})_{\nu})
               \sub(B_{-\nu-\mu})\sigma_0((B^{m-j})_{\mu})\\
      &=\int_{S^*M} dxd\xi\,
         \sum_{j=1}^{m-1}\,\frac1j\sum_{\nu}(-(\nu)_-) 
       \sigma_0((B^{j})_{\nu})
               \sub((B^{m-j})_{-\nu})     
\end{aligned}
$$%\end{equation}
which in view of Lemma \ref{lemsymm} equals
\begin{equation}
\label{CH1_s6sube(3)f}
%\begin{aligned}
    \int_{S^*M}  dxd\xi\,\sum_{\nu}(-(\nu_1)_-)
        \sum_{j=1}^{m-1}\frac1j\,
      \sub\big((B^{m-j})_{\nu}\big)\, 
              \sigma_0\big((B^j)_{-\nu}\big)
%\end{aligned}
\end{equation}
We make now in \eqref{CH1_s6sube(3)f} a change of index $j\ra{}m-j$
and get
\begin{equation}
\label{CH1_s6sube(3)ff}
    \int_{S^*M}  dxd\xi\,\sum_{\nu}(-(\nu_1)_-)
        \sum_{j=1}^{m-1}\frac1{m-j}\,
      \sub\big((B^{j})_{\nu}\big)\, 
              \sigma_0\big((B^{m-j})_{-\nu}\big)
\end{equation}
We sum \eqref{CH1_s6sube(3)ff} with \eqref{eq_100},
use $\frac1j+\frac1{m-j}=\frac{m}{j(m-j)}$,
refer to Lemma \ref{lemsymm}, and replace $j\ra{}m-j$ again
to conclude that \eqref{CH1_s6sube3} equals
\begin{equation}
\label{CH1_s6sube(3)fff}
     - m\int_{S^*M} \sum_{\nu}\nu_+
       \sum_{j=1}^{m-1}\,\,
      \, \sigma_0\bigg(\Big(\frac{B^{j}}{j}\Big)_{\nu}\bigg)\,
        \sub\bigg(\Big(\frac{B^{m-j}}{m-j}\Big)_{-\nu}\bigg)\,
           dxd\xi.
\end{equation}
It remains to notice that \eqref{CH1_s6sube(3)fff} in view of \eqref{eqEsub}
equals
$$%\begin{equation}
%\label{CH1_s6sube(3)ff4}
     \int_{S^*M}dxd\xi\,\sum_{k=1}^\infty k
         \,\int_0^{2\pi}\int_0^{2\pi} e^{ik(t_1-t_2)} \, 
   \tilde{W}_2[z^m](b_0^{t_1},b_0^{t_2})\,
               b_{{\rm sub}}^{t_2}\,\frac{dt_1}{2\pi}\frac{dt_2}{2\pi}.
$$%\end{equation}
\section
%[Contribution of the non-symmetric part of the subprincipal symbol]
{Contribution of the non-symmetric part of 
$\sub\!\big(B_{\k_1}\cdots{}B_{\k_m}\big)$}
%the subprincipal symbol of $B_{\k_1}\cdots{}B_{\k_m}$}
%{An auxiliary calculation for the third order Szeg\"o functional
%corresponding to $\UpsilonPoi$ 
%(non-symmetric contribution from the
%subprincipal symbol, Poisson brackets)}
\label{CH1_s6Poi}
Let us choose any $B\in\Psi^0(M)$.
The goal of this section is to calculate 
for an arbitrary $m\in\N$, $m\geq2$,
the expression
\begin{equation}
\label{CH1_s9Poi.e5Poi}
%\begin{aligned}
   \sum_{\k_1+\cdots+\k_m=0}\,
            M_m({\overline{\k}}) \,
   \int_{S^*M} dx\,d\xi \,\sum_{1\leq{}k<l\leq{}m}
                               \{\sigma_0(B_{\k_k}),\sigma_0(B_{\k_l})\}\,
      \prod_{\genfrac{}{}{0pt}{}{p=1}{p\neq{}k,p\neq{}l}}^m\sigma_0(B_{\k_p}).
%\end{aligned}
\end{equation}
The domain of summation is symmetric
in $\overline{\k}$. 
However the integrated over $S^*M$ sum is generally speaking not symmetric.
Each of the terms in the sum in \eqref{CH1_s9Poi.e5Poi}
possesses a partial symmetry.

Depending on $k$ and $l$ in \eqref{CH1_s9Poi.e5Poi}
there are three possible cases. In the first case 
all the indices except for those in the Poisson bracket
form {\em one\ }continuous block. In this case
it is convenient to rewrite $M_m({\overline{\k}})$ as follows.
Note that for $a,b\in\R$
\begin{equation}
\label{CH1_s6Poi.e1}
         \min(a,b)=a-(b-a)_-,
\end{equation}
let $p\in\N$, $p\geq2$, $\mu_1,\cdots,\mu_p\in\Z$
and set $s_j:=\mu_1+\cdots+\mu_j$, $j=1,\cdots,p$.
Then for any $j=1,\cdots,p-1$
$$%\begin{equation}
%\label{CH1_s6Poi.e2}
\begin{aligned}
   \null&\min(0,s_1,\cdots,s_j,s_{j+1},\cdots,s_p)\\ 
    &\quad=\min\big(\,\min(0,s_1,\cdots,s_j), 
    \min(s_j,s_{j+1},\cdots,s_m)\,\big)\\
    &\quad=\min\big(\,\min(0,s_1,\cdots,s_j), 
    s_j+\min(0,s_{j+1}-s_j,\cdots,s_p-s_j)\,\big).
\end{aligned}
$$%\end{equation}
Therefore it follows from \eqref{CH1_s6Poi.e1}
\begin{equation}
\label{eqONE}
\begin{aligned}  
     M_p(&\mu_{1},\cdots,\mu_{p}) 
      = M_{j}(\mu_{1},\cdots,\mu_{j})\\
  &- \big(\,\mu_{1}+\cdots+\mu_{j}
         -M_{j}(\mu_{1},\cdots,\mu_{j})
  +M_{p-j}(\mu_{{j+1}},\cdots,\mu_{p})\,\big)_-.
\end{aligned}
\end{equation}
The three subcases of the first
case are
$$
   k=1,\,l=2;\qquad k=1,\,l=m;\qquad k=m-1,\,l=m.
$$
In these three subcases using \eqref{eqONE} we rewrite 
$M_m(\k_{1},\cdots,\k_{m})$, respectively, as
$$
     M_{2}(\k_1,\k_2) - \big(\,\k_1+\k_2-M_{2}(\k_1,\k_2) 
           + M_{m-2}(\k_3,\cdots,\k_m)\,\big)_-,
$$
and
$$
\begin{aligned}
     -(\k_1)_- &- \Big(\,(\k_1)_+ + M_{m-2}(\k_2,\cdots,\k_{m-1})\\
                  &- \big(\, \k_2+\cdots+\k_{m-1}
        - M_{m-2}(\k_2,\cdots,\k_{m-1}) + \k_m\,\big)_-\,\Big)_-,
\end{aligned}
$$
and
$$
\begin{aligned}
     M_{m-2}(\k_1,\cdots,&\k_{m-2}) 
      - \big(\,\k_1+\cdots+\k_{m-2} \\
     &- M_{m-2}(\k_1,\cdots,\k_{m-2}) 
         + M_{2}(\k_{m-1},\k_m) \,\big)_-.
\end{aligned}
$$

In the second case all the indices but the two from the Poisson
bracket form {\em two\ }continuous blocks. 
In this case we rewrite $M_m({\overline{\k}})$ as follows.
Note that for $p\in\N$, $p\geq3$, $\mu_1,\cdots,\mu_p$ and
$s_1,\cdots,s_p$ as above, and any $j=1,\cdots,p-2$
$$%\begin{equation}
%\label{CH1_s6Poi.e2b}
\begin{aligned}
   \null&\min(0,s_1,\cdots,s_j,s_{j+1},\cdots,s_p)\\ 
    &\quad=\min\big(\,\min(0,s_1,\cdots,s_j), 
    \min(s_{j+1},\cdots,s_p)\,\big)\\
    &\quad=\min\big(\,\min(0,s_1,\cdots,s_j), 
    s_{j+1}+\min(0,s_{j+1}-s_j,\cdots,s_p-s_j)\,\big)
\end{aligned}
$$%\end{equation}
which together with \eqref{CH1_s6Poi.e1} implies
\begin{equation}
\label{eqTWO}
\begin{aligned}  
     \null&M_p(\mu_{1},\cdots,\mu_{p}) 
      = M_{j}(\k_{1},\cdots,\mu_{j})\\
 &- \big(\,\mu_{1}+\cdots+\mu_{j}
         -M_{j}(\mu_{1},\cdots,\mu_{j})
 +\mu_{{j+1}}
            +M_{p-j-1}(\mu_{{j+2}},\cdots,\mu_{p})\,\big)_-.
\end{aligned}
\end{equation}
There are three subcases:
first $k=1,l=r+1$, $2\leq{}r\leq{}m-2$, in which $M_m({\overline{\k}})$ 
after first using \eqref{eqONE} and then \eqref{eqTWO}
becomes
$$
\begin{aligned}
     -(\k_1)_- &- \Big(\,(\k_1)_+ + M_{r-1}(\k_2,\cdots,\k_{r})
- \big(\, \k_2+\cdots+\k_r\\
            &-M_{r-1}(\k_2,\cdots,\k_{r}) + \k_{r+1}
  + M_{m-r-1}(\k_{r+2},\cdots,\k_{m}) \,\big)_-\,\Big)_-,
\end{aligned}
$$
second: $k=r+1,l=r+2$, $1\leq{}r\leq{}m-3$, in 
which $M_m({\overline{\k}})$
 after first using \eqref{eqTWO} and then \eqref{eqONE}
becomes
$$
\begin{aligned}
      M_{r}&(\k_1,\cdots,\k_r) - \Big(\,
     \k_1+\cdots+\k_r-M_{r}(\k_1,\cdots,\k_r) + \k_{r+1}\\
    &-(\k_{r+2})_- - \big(\,
       (\k_{r+2})_+ + M_{m-r-2}(\k_{r+3},\cdots,\k_m)  \,\big)_-\,\Big)_-,
\end{aligned}
$$
third: $k=r+1,l=m$, $1\leq{}r\leq{}m-3$, in which $M_m({\overline{\k}})$ 
after first using first \eqref{eqTWO} and then \eqref{eqONE}
becomes
$$
\begin{aligned}
      M_{r}&(\k_1,\cdots,\k_r) - \Big(\,
     \k_1+\cdots+\k_r-M_{r}(\k_1,\cdots,\k_r) + \k_{r+1}\\
    &+ M_{m-r-2}(\k_{r+2},\cdots,\k_{m-1}) - \big(\,
       \k_{r+2}+\cdots+\k_{m-1} \\
    &-  M_{m-r-2}(\k_{r+2},\cdots,\k_{m-1}) + \k_m \,\big)_-\,\Big)_-.
\end{aligned}
$$

Finally in the third case the indices from the Poisson bracket are
taken to be $\k_{r+1}$ and $\k_{s+1}$,
where
$$
  1\leq r,\quad r+2\leq s,\quad s+2\leq m,
$$
and there are {\em three\ }continuous blocks
$1,\cdots,r$ and $r+2,\cdots,s$ and $s+2,\cdots,m$.
In that case $M_m({\overline{\k}})$ after using \eqref{eqTWO} twice
becomes
\begin{equation}
\label{CH1_s9Poi.e55}
\begin{aligned}
      M_{r}&(\k_1,\cdots,\k_r) - \Big(\,
     \k_1+\cdots+\k_r-M_{r}(\k_1,\cdots,\k_r) + \k_{r+1}\\
    &+ M_{s-r-1}(\k_{r+2},\cdots,\k_s) - \big(\,
       \k_{r+2}+\cdots+\k_s \\
    &-M_{s-r-1}(\k_{r+2},\cdots,\k_s) +\k_{s+1}
    + M_{m-s-1}(\k_{s+2},\cdots,\k_m)  \,\big)_-\,\Big)_-.
\end{aligned}
\end{equation}

We make the computation for the third case, the first and the second
are treated in the same way. We use the
a convenient reformulation
of the original form of the BSt (Theorem \ref{CH1_s1COR2}),
and one property of the $j$-maps $\Phi_j$, $j\in\N$, 
see Lemma \ref{CH1_s6Poi.l1} below. 
We rewrite the corresponding to the third
case part of \eqref{CH1_s9Poi.e5Poi} as
\begin{equation}
\label{CH1_s9Poi.e6}
\begin{aligned}
   \sum_{\k_1+\cdots+\k_m=0}\,
            &M_m({\overline{\k}}) \,
   \int_{S^*M} dx\,d\xi \,\sum_{r=1}^{m-3}\sum_{s=r+2}^{m-2}
                               \{\sigma_0(B_{\k_{r+1}}),\sigma_0(B_{\k_{s+1}})\}\\
           &\times\sigma_0(B_{\k_1}\cdots{}B_{\k_r})\,\sigma_0(B_{\k_{r+2}}\cdots{}B_{\k_s})
         \,\sigma_0(B_{\k_{s+2}}\cdots{}B_{\k_{m}}),
\end{aligned}
\end{equation}
where $k=r+1$, $l=s+1$, and each of the three products under
the $\sigma_0$ sign contains at least one factor (therefore this expression
is non-zero only for $m\geq5$). Now we observe that both 
domain of summation and each of the three products are symmetric
if we interchange the indices $\k_1,\cdots\k_r$,
$\k_{r+2},\cdots,\k_s$ and $\k_{s+2},\cdots,\k_m$ {\em separately,\ }that 
is within each of these three sets. Then their sums do not change,
and we can therefore consider the representation \eqref{CH1_s9Poi.e55}
as a function of 
$$
       M_{r}(\k_1,\cdots,\k_r),\quad 
      M_{s-r-1}(\k_{r+2},\cdots,\k_s),\quad
        M_{m-s-1}(\k_{s+2},\cdots,\k_m)      
$$
only. We interchange the indices within the three groups, use Theorem
\ref{CH1_s1COR2} and take in account 
Remark \ref{CH1_s4nnREM1}.
After that because of 
the symmetricity of domain of summation, and of the principal
symbols of the three products we may, and will, 
drop the coefficients $1/(r!)$,
$1/((s-r-1)!)$ and $1/((m-s-1)!)$ and the summations over the
permutations over the three groups of indices
(just as in Section \ref{CH1_s6aux}). 
We conclude that \eqref{CH1_s9Poi.e6} becomes
%\begin{equation}
%%\label{CH1_s9Poi.e65}
%\begin{aligned}
\begin{align}
\label{CH1_s9Poi.e65}
   \null&\sum_{\k_1+\cdots+\k_m=0}
        \sum_{r=1}^{m-3}\sum_{s=r+2}^{m-2}
         \sum_{\alpha=1}^r\frac1{\alpha!}
       \sum_{ \genfrac{}{}{0pt}{}{j_1,\cdots,j_\alpha\geq{}1}
                         {j_1+\cdots+j_\alpha=r} }
       \frac1{j_1\cdots{}j_\alpha}
 \sum_{\beta=1}^{s-r-1}\frac1{\beta!}
       \sum_{ \genfrac{}{}{0pt}{}{k_1,\cdots,k_\beta\geq{}1}
                 {k_1+\cdots+k_\beta=s-r-1} }
       \frac1{k_1\cdots{}k_\beta}\\
        &\times \sum_{\gamma=1}^{m-s-1}\frac1{\gamma!}
       \sum_{ \genfrac{}{}{0pt}{}{l_1,\cdots,l_\gamma\geq{}1}
        {l_1+\cdots+l_\gamma=m-s-1} }
       \frac1{l_1\cdots{}l_\gamma}
  \int_{S^*M} dx\,d\xi \,\,
    \big\{\sigma_0(B_{\k_{r+1}}),\sigma_0(B_{\k_{s+1}})\big\}\notag\displaybreak[0]\\
  &\times\Bigg[-(j_1(\overline{\k}))_- -\cdots-(j_\alpha(\overline{\k}))_- 
- \Big(\,
       (j_1(\overline{\k}))_+ +\cdots+(j_\alpha(\overline{\k}))_+ + \k_{r+1}\notag\displaybreak[0]\\
    & -(k_1(\overline{\k}))_- -\cdots-(k_\beta(\overline{\k}))_- - \big(\,
       (k_1(\overline{\k}))_+ +\cdots+(k_\beta(\overline{\k}))_+ +\k_{s+1}\notag\displaybreak[0]\\
  &-(l_1(\overline{\k}))_- -\cdots-(l_\gamma(\overline{\k}))_- \,\big)_-\,\Big)_-\Bigg]\notag\displaybreak[0]\\
    &\times\sigma_0(B_{\k_1}\cdots{}B_{\k_r})\,\sigma_0(B_{\k_{r+2}}\cdots{}B_{\k_s})
         \,\sigma_0(B_{\k_{s+2}}\cdots{}B_{\k_{m}})\notag.
\end{align}
%\end{aligned}
%\end{equation}
Here for all possible values of the indices {\em each\ }of the 
summation variables $\k_1,\cdots,\k_m$ is involved in some
of the $j_1(\overline{\k}),\cdots$, but only {\em once.\ }
We rewrite now the sum over $r$ and $s$ as a sum over three
summation variables $a,b,c\geq1$, which are the lengths of the
three continuous blocks, with the condition $a+b+c=m-2$.
Note also that it does not matter anymore on which place 
$\k_{r+1}$ and $\k_{s+1}$ stands, it is only important 
that the same letter is
used in the non-positive valued function coming from the
square bracket in \eqref{CH1_s9Poi.e65}. 
Let us make two changes of variables
$$
  \k_1:= \k_{r+1},\quad \k_2:= \k_{r+1}
$$
and for all $a,b,c\geq1$ with $a+b+c=m-2$ 
$$
\begin{aligned}
       \mu_1&:=j_1(\overline{\k}),\cdots,\mu_\alpha:=j_\alpha(\overline{\k}),
                      \qquad 1\leq\alpha\leq a\\
       \nu_1&:=k_1(\overline{\k}),\cdots,\nu_\beta:=k_\beta(\overline{\k}),
                      \qquad 1\leq\beta\leq b\\
       \rho_1&:=l_1(\overline{\k}),\cdots,\rho_\gamma:=l_\gamma(\overline{\k}),
                      \qquad 1\leq\gamma\leq c.
\end{aligned}
$$
Then the
square bracket in \eqref{CH1_s9Poi.e65} becomes exactly 
the defined by \eqref{CH1_s5ieOj3}
non-positive valued function $\O_{\alpha,\beta,\gamma}^{(3)}$.
Now we carry out the $((a-\alpha)+(b-\beta)+(c-\gamma))$ free
summations over the rest of $\k$'s, just as in Section \ref{CH1_s6aux}.
We get
% \begin{equation}
%%\label{CH1_s9Poi.e68}
\begin{align}
%\label{CH1_s9Poi.e68}
   \null&\int_{S^*M}dxd\xi \,\sum_{ \genfrac{}{}{0pt}{}{a,b,c\geq1}{a+b+c=m-2} }
       \sum_{\alpha=1}^a\sum_{\beta=1}^b\sum_{\gamma=1}^c
          \frac1{\alpha!}\frac1{\beta!}\frac1{\gamma!}\displaybreak[0]\\
   &\times\sum_{\genfrac{}{}{0pt}{}{\k_1+\k_2+\mu_1+\cdots+\mu_\alpha+}
       {+\nu_1+\cdots+\nu_\beta+\rho_1+\cdots+\rho_\gamma=0}}
   \Omega_{\alpha,\beta,\gamma}^{(3)}
         \big(\k_1,\k_2,\mu_1,\cdots,\mu_\alpha,\nu_1,\cdots,\nu_\beta,
        \rho_1,\cdots,\rho_\gamma\big)
       \notag\displaybreak[0]\\
   &\times\,\big\{\sigma_0(B_{\k_1}),\sigma_0(B_{\k_2})\big\}\,\sigma_0\Bigg(       
      \sum_{ \genfrac{}{}{0pt}{}{j_1,\cdots,j_\alpha\geq{}1}
                       {j_1+\cdots+j_\alpha=a} }
       \frac{\big(B^{j_1}\big)_{\mu_1}\cdots
                       \big(B^{j_\alpha}\big)_{\mu_\alpha}}
              {j_1\cdots{}j_\alpha}\Bigg)\notag\displaybreak[0]\\
   &\times\sigma_0\Bigg(
       \sum_{ \genfrac{}{}{0pt}{}{k_1,\cdots,k_\beta\geq{}1}
                     {k_1+\cdots+k_\beta=b} }
       \frac{\big(B^{k_1}\big)_{\nu_1}\cdots
                       \big(B^{k_\beta}\big)_{\nu_\beta}}
         {k_1\cdots{}k_\beta}\Bigg)
\,\sigma_0\Bigg(
       \sum_{ \genfrac{}{}{0pt}{}{l_1,\cdots,l_\gamma\geq{}1}
        {l_1+\cdots+l_\gamma=c} }
       \frac{\big(B^{l_1}\big)_{\rho_1}\cdots
                       \big(B^{l_\gamma}\big)_{\rho_\gamma}}
               {l_1\cdots{}l_\gamma}\Bigg).\notag
\end{align}
%\end{aligned}
%\end{equation}
Now by the definition of the Fourier coefficient,
by Egorov's theorem and the fact that $\sigma_0(B^j)=(\sigma_0(B))^j$, $j\in\N$,
we obtain, just as in Section \ref{CH1_s6aux},
%(now we write $j,k,l$ instead of $\alpha,\beta,\gamma$),
\begin{equation}
\label{CH1_s6Poi.e69}
\begin{aligned}
   \null&\int_{S^*M}dxd\xi 
         \,\sum_{\genfrac{}{}{0pt}{}{a,b,c\geq1}{a+b+c=m-2}}
       \sum_{\alpha=1}^a\sum_{\beta=1}^b\sum_{\gamma=1}^c
          \frac1{\alpha!}\frac1{\beta!}\frac1{\gamma!}\\
   &\times\sum_{\genfrac{}{}{0pt}{}{\k_1+\k_2+\mu_1+\cdots+\mu_\alpha+}
          {+\nu_1+\cdots+\nu_\beta+\rho_1+\cdots+\rho_\gamma=0}}
   \Omega_{\alpha,\beta,\gamma}^{(3)}
         \big(\k_1,\k_2,\mu_1,\cdots,\mu_\alpha,\nu_1,\cdots,\nu_\beta,
        \rho_1,\cdots,\rho_\gamma\big)
       \\
   &\times\,\big\{\sigma_0(B_{\k_1}),\sigma_0(B_{\k_2})\big\}\\ 
&\times\int_0^{2\pi}\!\cdots\!\int_0^{2\pi}
\,\frac{dr_1}{2\pi}\cdots\frac{dr_\alpha}{2\pi}
\,\frac{ds_1}{2\pi}\cdots\frac{ds_\beta}{2\pi}
\,\frac{dt_1}{2\pi}\cdots\frac{dt_\gamma}{2\pi}\\
&\times e^{i(\mu_1r_1+\cdots+\mu_\alpha{}r_\alpha
  +\nu_1s_1+\cdots+\nu_\beta{}s_\beta
       +\rho_1t_1+\cdots+\rho_\gamma{}t_\gamma)}\\
&\times
       \Phi_\alpha[z^a](b_0^{r_1},\cdots,b_0^{r_\alpha})\,
       \Phi_\beta[z^b](b_0^{s_1},\cdots,b_0^{s_\beta})\,
        \Phi_\gamma[z^c](,b_0^{t_1},\cdots,b_0^{t_\gamma}).
\end{aligned}
\end{equation}
We use the fact that $\Phi_\alpha[z^a]=0$ for $\alpha>a$
to extend the summation over each of the variables
$\alpha,\beta,\gamma$ to the whole of $\N$, and make these
summations the outer ones. Finally we use the following important
property
of the $j$-maps $\Phi_j$, $j\in\N$,
whose proof immediately follows from the definition of $\Phi_j$
and is omitted.
\begin{lem}
\label{CH1_s6Poi.l1} 
Let choose an arbitrary $p\in\N$.
For any fixed real or complex $x_j$, $y_k$, $z_l$, 
$j,k,l=1,\cdots,p-2$ the following identity holds
\begin{equation}
\label{s6Poie0000}
\begin{aligned}
  \sum_{\genfrac{}{}{0pt}{}{a,b,c\geq1}{a+b+c=p}}
       &\sum_{\alpha=1}^a\sum_{\beta=1}^b\sum_{\gamma=1}^c
          \frac1{\alpha!}\frac1{\beta!}\frac1{\gamma!}\\
   &\quad\times\Phi_\alpha[z^a](x_1,\cdots,x_\alpha)\,
       \Phi_\beta[z^b](y_1,\cdots,y_{\beta})\,
        \Phi_\gamma[z^c](z_1,\cdots,z_{\gamma})\\
   =\sum_{\alpha=1}^\infinity
          &\sum_{\beta=1}^\infinity
 \sum_{\gamma=1}^\infinity
          \frac1{\alpha!}\frac1{\beta!}\frac1{\gamma!}\\
     &\qquad\times\Phi_{\alpha+\beta+\gamma}
         [z^{p}](x_1,\cdots,x_\alpha,y_1,\cdots,y_{\beta},
                z_1,\cdots,z_{\gamma}).
\end{aligned}
\end{equation}
\end{lem}
Note that the sums terminate for $\alpha+\beta+\gamma>p$,
and $\alpha,\beta,\gamma\geq1$. Therefore 
only $p$ of the variables $x_1,\cdots,z_{p-2}$ 
are present at the each term on the right-hand side in \eqref{s6Poie0000}.

The formula \eqref{CH1_s5ie1Phi20Poi3} now 
follows from \eqref{CH1_s6Poi.e69} and Lemma
\ref{CH1_s6Poi.l1} with $p=m-2$. 
It is important that we have reduced the number
of $j$-maps from three to one, and that the $j$-map $\Phi_j[f]$ 
is linear in $f$ for all $j\in\N$. Therefore \eqref{CH1_s5ie1Phi20Poi3}
holds for any $f\in\an$. We divide by $z^2$, because
for a monomial $f(z)=z^m$, $m\geq5$, 
the $j$-map should be evaluated at $z^{m-2}$.
We subtract the fourth degree Taylor polynomial $T_4[f]$ of $f$ about $t=0$ 
because
the term under consideration is absent for all polynomials
$f$ of degree $\leq4$. 

The formulas \eqref{CH1_s5ie1Phi20Poi1} and \eqref{CH1_s5ie1Phi20Poi2} 
are proved analogously.
Because the corresponding expressions appear for the monomials
of degree at least $3$ and $4$, respectively,
we subtract $T_2[f]$ and $T_3[f]$, respectively.
\section{Proof of Proposition~\ref{prop3}}
\label{secpfprop3}
We will need the following statement (see for instance \cite{GR}).
\begin{lem}
\label{lem_z}
Let $n,m\in\N$. Then the following holds, as $n\ra\infty$,
$$
\begin{aligned}
  \sum_{k=1}^n k^m 
           &= \frac{n^{m+1}}{m+1}+\frac12n^{m}
           +\frac{m}{12}n^{m-1}+ O(n^{m-2})\\
  \sum_{k=1}^n k^{-1} &= \log n +\gamma+\frac12n^{-1}+ O(n^{-2})\\
  \sum_{k=1}^n k^{-2} &= \zeta(2) -\frac1n+ O(n^{-2})\\
  \sum_{k=1}^n k^{-m} &= \zeta(m) + O(n^{-m+1}),\qquad m\geq3.
\end{aligned}
$$
\end{lem}
Note that 
\begin{equation}
\label{eqstrace}
   \tr(P_nG)=\sum_{k=1}^n\tr(\pi_kG),\qquad n\in\N.
\end{equation}
Assume first that $d\geq4$. Then by \eqref{eq_CdV}
$$%\begin{equation}
%\label{eq_3CdV}
         \Big| \tr(\pi_kG)
   - k^{d-1}\,R_0(G) - k^{d-2}\,R_1(G) - k^{d-3}\,R_2(G) \Big|
\leq c_1(G)\,k^{d-4}, \quad k\in\N,
$$%\end{equation}
and so in view of \eqref{eqstrace}, for any $n\in\N$,
\begin{equation}
\label{eq_4CdV}
\begin{aligned}
 \null\Big| \tr(P_nG) &- n^d\cdot\frac1d\,R_0(G) -
  n^{d-1}\cdot\Big(\frac12\,R_0(G) + \frac1{d-1}\,R_1(G)\Big)\\
&- 
  n^{d-2}\cdot\Big(\frac{d-1}{12}\,R_0(G) + \frac12\,R_1(G)
            + \frac1{d-2}\,R_2(G)\Big)\\
%+ \log n\cdot R_d(G) \\
&+ c_3(R_0(G),R_1(G),R_2(G),d)\cdot n^{d-3}\Big| 
          \leq c_2(G,d)\,n^{d-3}.
\end{aligned}
\end{equation}
For $d=3$ there also appears a term with $\log n$
in the left-hand side of \eqref{eq_4CdV}.
This proves part (iii).

For $d=1,2$ when we sum over $k=1,\cdots,n$ in \eqref{eq_CdV}
there is a subtle point, namely, 
the constant coefficient in \eqref{eqstrace}, as $n\ra\infty$. 
The terms of all orders
in \eqref{eq_CdV}, and also the possible rapidly decaying term, 
will contribute to it. 
Assume \eqref{eq_assu}. Then the following series is absolutely concergent
%definition makes sense 
\begin{equation}
\label{eq_eps}
   \eps_k(G):= \tr(\pi_kG)
   -\sum_{l=0}^{+\infty} k^{d-1-l}\,R_l(G),\qquad k\in\N.
\end{equation}
Furthermore, for any $N\in\N$ and all $k\in\N$ by \eqref{eq_CdV}
\begin{equation}
\label{eq_epseps}
\begin{aligned}
   |\eps_k(G) | \leq \big| \tr(\pi_kG)
   &-\sum_{l=0}^{N+d-2} k^{d-1-l}\,R_l(G)\big| 
       +\big| \sum_{l=N+d-1}^{\infty} k^{d-1-l}\,R_l(G)\big| \\
   &\leq c_N(G)\,k^{-N} + k^{-N}\sum_{l=0}^\infty|R_l(G)| 
       \leq \tilde{c}_N(G)\cdot k^{-N}
\end{aligned}
\end{equation}
in view of \eqref{eq_assu}.
Note that $C(G)$ defined by \eqref{eq_C} equals
\begin{equation}
\label{eq_Ceps}
   C(G)= \sum_{k=1}^\infty \eps_k(G),
\end{equation}
the series being absolutely convergent.
Now summing over $k=1,\cdots,n$ in \eqref{eq_eps}
we obtain %and interchanging the finite and the infinite sum
\begin{equation}
\label{eq_eps1}
   \tr(P_nG) = 
   \sum_{l=0}^{+\infty} R_l(G) \sum_{k=1}^n k^{d-1-l}
          + \sum_{k=1}^n\eps_k(G).
\end{equation}
Because $\eps_k(G)$ decays rapidly \eqref{eq_epseps}, as $k\ra\infty$,
and by \eqref{eq_Ceps}, $\sum_{k=1}^n\eps_k(G)$ 
converges rapidly to $C(G)$. 
From this and \eqref{eq_eps1} we can obtain an asymptotics
of $\tr(P_nG)$ up to any negative order. In particular,
taking into account Lemma~\ref{lem_z} we prove (i) and (ii).
\section%[Invariant formulas for $\Phi_j$ and $F_{j}$, $j\in\N$]
{A formula for $\Phi_j[f]$, $f\in\an$, $j\in\N$}
%Construction of the invariant formulas for
%the functors $\tilde\Phi_j$, $\Phi_j$ and $F_{j}$, $j\in\N$}
\label{CH1_s6Phi}
We find a formula for $\Phi_j[f]$, $f\in\an$, 
in terms of an auxiliary linear $j$-map 
$\tilde\Phi_j$, $j\in\N$, see \eqref{CH1_s6Phi.e2p}.
The latter acts on monomials as follows 
\begin{equation}
\label{CH1_s6Phi.e1}
           \tilde\Phi_j[z^m](\xi_1,\cdots,\xi_j) :=
                 \sum_{\genfrac{}{}{0pt}{}{l_1,\cdots,l_j\geq0}
                                    {l_1+\cdots+l_j=m}} 
                         \xi_1^{l_1}\cdots\xi_j^{l_j},\quad m=0,1,2,\cdots,
\end{equation}
which is the complete symmetric function of degree $m$
evaluated at the point $(\xi_1,\cdots,\xi_j,0,\cdots)$.
To write a formula for $\tilde\Phi_j[f]$, $f\in\an$,
we use an idea suggested by Kurt Johansson.
It uses a Cauchy integral representation 
of \eqref{CH1_s6Phi.e1} 
via the generating function. Namely, by the identity \cite[(I.2.5)]{M}
for $(\xi_1,\cdots,\xi_j,0,\cdots)$
$$
   \sum_{n=0}^\infty\tilde\Phi_j[z^n](\xi_1,\cdots,\xi_j)\cdot\frac{1}{\zeta^n}
       = \prod_{k=1}^j\frac{\zeta}{\zeta-\xi_k}.
$$
After a multiplication by $\zeta^{m-1}$ and an integration over a 
contour $\gamma\subset\C$ which circumferences 
the points $0,\xi_1,\cdots,\xi_j$ we single out $\tilde\Phi_j[z^m]$
and obtain
$$%\begin{equation}
%\label{CH1_s6Phi100}
   \tilde\Phi_j[z^m](\xi_1,\cdots,\xi_j) = \frac1{2\pi{}i}\int_\gamma
                 \frac{\zeta^{j-1}}{(\zeta-\xi_1)\cdots(\zeta-\xi_j)}\,\zeta^m\,d\zeta.
$$%\end{equation}
Therefore 
for $f\in\A_0$ (an analytic on $\C$ function which might have a constant
term) we can define
$$%\begin{equation}
%\label{CH1_s6Phi99}
   \tilde\Phi_j[f(z)](\xi_1,\cdots,\xi_j) = \frac1{2\pi{}i}\int_\gamma
\frac{\zeta^{j-1}}{(\zeta-\xi_1)\cdots(\zeta-\xi_j)}\,f(\zeta)\,d\zeta.
$$%\end{equation}

The $j$-map $\Phi_j$ is now defined for any $f\in\an$ by 
\begin{equation}
\label{CH1_s6Phi.e2p}
       \Phi_j[f(z)](x_1,\cdots,x_j) = \int_0^{x_1}\cdots\,\int_0^{x_j}
              \tilde\Phi_j[z^{-j}f(z)](\xi_1,\cdots,\xi_j)\,d\xi_1\cdots d\xi_j.
\end{equation}
Now we take \eqref{CH1_s6Phi.e2p} as a definition of $\Phi_j$,
then \eqref{CH1_s5iePhi2} with $m\geq{}j$ holds, and it only
remains to prove that $\Phi_j$
vanishes on the set of polynomials of degree $j-1$
with no constant term. 
Note that the integrand in \eqref{CH1_s6Phi.e2p} equals
\begin{equation}
\label{CH1_s6Phi100}
     \frac1{2\pi{}i}\int_\gamma
        \frac{\zeta^{-1}f(\zeta)}{(\zeta-\xi_1)\cdots(\zeta-\xi_j)}\,d\zeta.
\end{equation}
The integral does not depend on the contour of integration,
if only all $\xi_1,\cdots$, $\xi_j$ are inside it.
Let $\gamma=\{\zeta:|\zeta|=R\}$, $R\ra\infty$.
If $f(z)=z^k$, $k=1,\cdots,j-1$, then the absolute value of \eqref{CH1_s6Phi100}
is estimated by $(2\pi)^{-1}R^{j-2}\cdot{}2\pi{}R/R^j=R^{-1}\ra0$,
as $R\ra\infty$. Therefore \eqref{CH1_s6Phi100} must be $0$.

If $f\in\an$ then $z^{-1}f(z)$ is analytic. Therefore \eqref{CH1_s6Phi100}
equals the sum of $j$ residues at the points $\xi_1,\cdots,\xi_j$.
For instance, for $j=2$ and any $f\in\an$
$$
   \Phi_2[f](x_1,x_2) = \int_0^{x_1}\int_0^{x_2}
   \frac{\xi_1^{-1}f(\xi_1)-\xi_2^{-1}f(\xi_2)}{\xi_1-\xi_2}\,d\xi_1 d\xi_1.
$$
%This justifies the definitions \eqref{CH1_s5ieFprin}.
\begin{rem}
Another way to construct
$\tilde\Phi_j$, $j\in\N$, is by induction on $j$.
In that case one uses a simple formula 
\begin{equation}
\label{CH1_s7e7}
  \frac{u^{r+1} - v^{r+1}}{u - v} = \sum_{ \genfrac{}{}{0pt}{}{p,q\geq0}{p+q=r} }
          u^pv^q, \qquad u,v\in\C,
\end{equation}
for $r=0,1,2,\cdots$. This derivation is however
longer than the above argument. 
The formula \eqref{CH1_s7e7} was used
in \cite{LRS} for the computation of the 2-map $W_2$ defined in
\eqref{CH1_s5ieW}. In that case the induction is not needed.
\end{rem}
Finally we find a formula for $W_3[f]$, $f\in\an$,
where the action $W_3[z^m]$, $m\in\N$, is given 
by \eqref{CH1_s6e16biss}.
In view of \eqref{CH1_s6Phi.e1} 
using the integration as in \eqref{CH1_s6Phi.e2p}
and moving out $x_3$ we obtain
$$
       W_3[f(z)](x_1,x_2,x_3) = x_3\int_0^{x_1}\int_0^{x_2}
           \tilde{\Phi}_3[z^{-3}f(z)](\xi_1,\xi_2,x_3)\,d\xi_1d\xi_2.
$$
Now evaluating the three residues in \eqref{CH1_s6Phi100}
for $j=3$
we obtain~\eqref{CH1_s5ieFprin}.
%
%\setcounter{equation}{0}
%\setcounter{defin}{0}
%\section[The formula gHD, and its derivation from BSt]
%{The formula gHD, and its derivation from BSt}
%{Generalized Hunt--Dyson formula, and its derivation
%from Bohnenblust--Spitzer theorem}
\section{Formulation of the gHD, and of a %a combinatorial formula CF
version of the BSt%), and of gHD
}
\label{CH1_s4nn}
In this section we state the gHD
and a convenient for our purposes
version of the BSt, see \cite{GLc} % and Remark~\ref{rem_comb}
for details. Let $m\in\N$ and 
$\k_1,\ldots,\k_m\in\R$.
Recall the notation \eqref{eq_Mm},
and  for each permutation $\tau\in{}S_m$ write
$$
    {\overline{\k}}_\tau  := (\k_{\tau_1},\ldots,\k_{\tau_m}).
$$
Fix any $j=1,\cdots,m$. Recall \cite{C} that a {\em partition\ }of $m$
is a way to write $m$ as a sum of natural numbers with no respect
to the order of the summands. These summands are called 
{\em parts,\ }and their values are called {\em lengths}. 
Note that each of the parts has length at least one.
For an arbitrary partition of  $m$ into $j$ parts
$$
           k_1\geq1,\cdots,k_j\geq1,\qquad k_1+\cdots+k_j=m,
$$
we introduce the notation
\begin{equation}
\label{CH1_s1e6prim}
\begin{aligned}
   k_1({\overline{\k}}_\tau) &:= \k_{\tau_1}+\cdots+\k_{\tau_{k_1}}\\
   k_2({\overline{\k}}_\tau) &:= \k_{\tau_{k_1+1}}+\cdots+
           \k_{\tau_{k_1+k_2}}\\
                 &\quad\cdots\\
   k_j({\overline{\k}}_\tau) &:= \k_{\tau_{k_1+\cdots+k_{j-1}+1}}+\cdots+
         \k_{\tau_{k_1+\cdots+k_{j-1}+k_j}}.
\end{aligned}
\end{equation}
Each of $k_l({\overline{\k}}_\tau)$, $l=1,\cdots,j$, is a sum
of $k_l$ permuted variables out of $\k_{\tau_1},\cdots,\k_{\tau_m}$ 
so that each of the permuted variables enters exactly one sum.
Note also that because $k_1+\cdots+k_j=m$ one has
$$
\begin{aligned}
       k_1({\overline{\k}}_\tau) +  \cdots + k_j({\overline{\k}}_\tau) 
             &= \k_{\tau_1}+\cdots+\k_{\tau_m}\\
        &= \k_{1}+\cdots+\k_{m}.
\end{aligned}
$$
For $a\in\R$, denote $-(a)_-:=\min(0,a)$ and $(a)_+:=\max(0,a)$.
We state now a combinatorial formula, called
{\em CF,\ }which is an equivalent
form of the BSt \cite[Theorem~2.2]{S1}, and does not involve any
advanced combinatorial coefficients. 
This formula is very suitable % well adopted 
for the calculation of sums which arise in Section~\ref{CH1_s6Poi} 
%higher order Szeg\"o terms
and in computations of the joint distributions for random walks 
as in Theorem \ref{CH1_s5ith1RW}. Recall the notation \eqref{eq_Mm}.
%, it is our main tool in further constructions.
\begin{thm}[CF: an equivalent version of the BSt]
\label{CH1_s1COR2}
For any $m\in\N$, arbitrary $\k_1,\cdots,\k_m\in\R$,
and any real- or complex-valued function $f$ defined on 
the left half-axis, the following holds
\begin{equation}
\label{CH1_s1eCOR2}
\begin{aligned}
  \sum_{\tau\in{}S_m}\,
         &f(\,M_m(\overline{\k}_\tau)\,)\\
          &=  \sum_{\tau\in{}S_m}\, \sum_{j=1}^m\frac1{j!}
       \sum_{ \genfrac{}{}{0pt}{}{k_1,\cdots,k_j\geq1}{k_1+\cdots+k_j=m} }
          \,  \frac{f\big(\,-(k_1(\k_\tau))_--\cdots-(k_j(\k_\tau))_-\,\big)}
                          {k_1\cdots k_j}
\end{aligned}
\end{equation}
\end{thm}
This holds because the {\em sets of the values of the arguments\ }of 
$f$ on the left- and
on the right-hand side in \eqref{CH1_s1eCOR2} contain the same
numbers with the same multiplicities, by the BSt,
the rest being just an account of the number 
of conjugacy classes in $S_m$,
see Section~3 and especially the proof of Lemma~3.3 in \cite{GLc}.
%The reader who wishes to recall what the BSt is, and why the CF 
%is an equivalent form of BSt, is referred to 
%and \cite{GLc} and \cite[Remark~2.1.4]{GiPhD}.

For a monomial $f(z)=z^n$, $n\in\N$,
a further calculation in CF can be carried out, of one subtracts 
the value of $f$ at the ``previous'' maximum. 
The right-hand side of the  resulting formula
has a {\em multiplicative,\ }and not
{\em additive,\ } as in the CF, form,
which is important for the calculation of the 
sums of convolution type in Section~\ref{CH1_s6aux}. 
\begin{thm}[Generalized Hunt--Dyson formula]
\label{CH1_S2TH1}
For any power $n\in\N$,
an arbitrary number of variables $m\in\N$
and for arbitrary fixed values of the real variables $\k_1,\cdots,\k_m$
one has
\begin{equation}
\label{CH1_s2e10}
\begin{aligned}
  \sum_{\tau\in{}S_m}\,
          &\ls \lp M_{m}({\overline{\k}}_\tau) \rp^n
                             - ( M_{m-1}({\overline{\k}}_\tau) )^n \rs \\
         &=  \sum_{\tau\in{}S_m}\, \sum_{j=1}^{\min(m,n)}
          \,\frac1{j!}\,
 \sum_{   \genfrac{}{}{0pt}{}  {  k_1,\cdots,k_j \geq 1 }
                  { k_1+\cdots+k_j = m }
           } \,
  \sum_{   \genfrac{}{}{0pt}{}  {  l_1,\cdots,l_j \geq 1 }
                  { l_1+\cdots+l_j = n }
           } \,
  \lp\genfrac{}{}{0pt}{} {n}{ l_1,\cdots, l_j} \rp \\
      &\qquad\times\,
               \frac{(-(k_1({\overline{\k}}_\tau))_-)^{l_1}}{k_1}\cdot\ldots\cdot
                \frac{(-(k_j({\overline{\k}}_\tau))_-)^{l_j}}{k_j}.
\end{aligned}
\end{equation}
\end{thm}
\begin{rem} In the case $n=1$ the identity \eqref{CH1_s2e10} 
becomes the usual 
Hunt--Dyson formula \cite[after~(4.8)]{K} 
\begin{equation}
\label{CH1_s2e10orig}
\begin{aligned}
  \sum_{\tau\in{}S_m}\,
          \ls \,M_m({\overline{\k}}_\tau) 
                             - M_{m-1}({\overline{\k}}_\tau)\,\rs %\\
         &=  \sum_{\tau\in{}S_m}\,
            \frac{-(\k_{\tau_1}+\cdots+\k_{\tau_m})_-}{m}\\
          &=(m-1)!\,(-(\k_{1}+\cdots+\k_{m})_-).
\end{aligned}
\end{equation}
\end{rem}
\begin{rem} 
\label{CH1_s4nnREM1}
It is important that Theorem \ref{CH1_s1COR2} and
\ref{CH1_S2TH1} holds if the symmetric group
$S_m$ is replaced with a {\em larger\ }symmetric group of indices.
One does however need
the summation over the {\em whole\ }group $S_m$. The usual
formula HD holds for the sums over the cyclic
subgroup $C_m$, as well. %, at the least. 
\end{rem}
\begin{rem}
The usual formula HD
of course can be
viewed as a corollary of the Bohnen\-blust--Spitzer identity.
However one needs the 
formula \eqref{CH1_s2e10} for {\em all\ }$n\in\N$ 
to reprove the BSt as is done in \cite[Section~4]{GLc}.
\end{rem}
\begin{rem}
The BSt, the gHD, and the HD
hold if the minima are replaced with the maxima,
and the negative parts $-(\cdot)_-$ are replaced with $(\cdot)_+$.
\end{rem}
%
%% End of article:

%% optional
% Appendices

% Appendix without title:
%\appendix{}

% Appendix with title:
%\appendix{Title}

% Appendix with letter:
%\appendix{B}

% Appendix with letter and title:
%\appendix{C}
%\appendixtitle{This is an appendix title}

%% optional

%\bibliographystyle{plain}
%\bibliography{atmpl}
%\bibitem[B1]{B1}
%\bibitem[W4]{W3}
%\end{article}
%\end{document}
%% Not optional, necessary:
%\begin{references}
\bibliographystyle{plain}

%
%\end{references}

%% This command is necessary! ==>>
%\end{article}
\end{document}